\documentclass[12pt,a4paper]{article}

\usepackage[margin=1in]{geometry}
\usepackage[T1]{fontenc}
\usepackage{lmodern}
\usepackage{setspace}
\usepackage[round,authoryear]{natbib}
\usepackage{amsmath,amssymb,amsthm,mathtools,bm}
\usepackage{booktabs}
\usepackage{enumitem}
\usepackage{hyperref}
\usepackage{url}
\usepackage{graphicx}
\usepackage{algorithm}
\usepackage{algpseudocode}
\usepackage{xcolor}
\usepackage{caption}
\hypersetup{hidelinks}

\newtheorem{theorem}{Theorem}
\newtheorem{proposition}{Proposition}
\newtheorem{lemma}{Lemma}
\newtheorem{corollary}{Corollary}
\newtheorem{definition}{Definition}
\newtheorem{remark}{Remark}

\newcommand{\E}{\mathbb E}
\newcommand{\Pp}{\mathbb P}
\newcommand{\R}{\mathbb R}

\newcommand{\Var}{\operatorname{Var}}
\newcommand{\Cov}{\operatorname{Cov}}

\newcommand{\calF}{\mathcal F}
\newcommand{\calG}{\mathcal G}
\newcommand{\calX}{\mathcal X}
\newcommand{\calA}{\mathcal A}

\newcommand{\op}{\operatorname{op}}
\newcommand{\asto}{\xrightarrow{\mathrm{a.s.}}}

\newcommand{\norm}[1]{\left\lVert #1\right\rVert}

\newcommand{\trans}{^{\mathsf T}}

\newcommand{\bracket}[1]{\left\langle #1\right\rangle}
\newcommand{\qv}[1]{\left[ #1\right]}

\title{Martingalized U-statistics-based sequential inference: a blessing of degeneracy}
\author{Weijia Li\thanks{Department of Statistics and Data Science, Tsinghua University.}\hspace{0.4em},
Leheng Cai\footnotemark[1]\hspace{0.6em}and
Qirui Hu\thanks{\hangindent=1.8em\hangafter=1 School of Statistics and Data Science, Shanghai University of Finance and Economics.\\
Department of Mathematics, Ruhr University Bochum}\hspace{0.4em}\textsuperscript{,}\hspace{0.1em}\thanks{Corresponding author: huqirui@mail.shufe.edu.cn}\hspace{0.4em}}
\date{}

\begin{document}
\maketitle

\begin{abstract}
Degenerate degree-two U-statistics have non-Gaussian limits, making sequential inference substantially more difficult than in the nondegenerate case. However, row normalization preserves the target parameter and makes the centered numerator a martingale on a harmonic variance clock, at the cost of logarithmic variance inflation. We establish a strong Gaussian approximation for this process and consistently estimate its variance through online updates, allowing spectrum-free calibration. Combining this strong approximation with Gaussian-mixture (GM) and stitched law-of-the-iterated-logarithm (LIL) boundaries yields confidence sequences with asymptotic time-uniform coverage and corresponding sequential tests for null hypotheses under which the kernel is degenerate. We further extend the framework to multi-dimensional vector-valued kernels, yielding Mahalanobis confidence ellipsoids based on consistent covariance estimation. Numerical experiments illustrate sequential type-I error control and favorable detection power in
two-sample and independence testing against spectral baselines, and in
three-sample homogeneity testing using vector confidence ellipsoids.

\end{abstract}

\section{Introduction}
\label{sec:introduction}

Many kernel- and distance-based tests are based on degree-two U- or
V-statistics that become degenerate under the null hypothesis. Prominent
examples include the maximum mean discrepancy (MMD)
\citep{GrettonEtAl2012}, the Hilbert--Schmidt independence criterion (HSIC)
\citep{GrettonEtAl2005}, energy statistics, and distance covariance
\citep{SzekelyRizzoBakirov2007,SzekelyRizzo2013,SejdinovicEtAl2013}.
For degree-two U-statistics, the Hoeffding decomposition yields
fundamentally different asymptotic behavior according to whether the
first-order projection vanishes. In the nondegenerate case, the first-order
term gives the usual Gaussian limit, whereas under degeneracy the normalized
statistic typically converges to a quadratic Gaussian chaos whose law
depends on the spectrum of the associated kernel operator
\citep{Hoeffding1948,Serfling1980,Neuhaus1977,
DehlingDenkerPhilipp1984,ArconesGine1993,ArconesGine1995,
KanagawaYoshihara1994}. The spectral dependence of the degenerate case complicates
fixed-sample calibration and has motivated, among other approaches,
permutation procedures and modifications of the original statistic that
restore Gaussian limits
\citep{KimRamdas2024,ShekharKimRamdas2022,ShekharKimRamdas2023}.

Sequential monitoring, in which inference may be updated after each new
observation, requires a stronger form of error control than fixed-sample
inference. Rather than controlling the statistic at a single deterministic
sample size, one must control the probability that it crosses a prescribed
boundary at any monitoring time. Classical sequential analysis studies such
boundary-crossing events through mixture constructions and related methods
\citep{Robbins1970a,Robbins1970b}, while modern anytime-valid inference is
often based on nonnegative supermartingales, e-processes, mixture boundaries,
stitching, and betting
\citep{HowardRamdasMcAuliffeSekhon2020,
HowardRamdasMcAuliffeSekhon2021,RamdasGrunwaldVovkShafer2023}.

A complementary asymptotic approach derives time-uniform inference from
process-level approximations. If a statistic can be coupled to a Gaussian
process with an error negligible relative to the relevant boundary scale,
Gaussian boundary-crossing results can be transferred to the original
process. This principle underlies the time-uniform central limit theory and
asymptotic confidence-sequence framework of
\citet{WaudbySmithArbourSinhaKennedyRamdas2024}; related sharp asymptotic
constructions are developed by \citet{GnettnerKirch2025}. For ordinary
degenerate U-statistics, however, the natural approximating process is a
spectrum-dependent quadratic Gaussian chaos rather than a Gaussian process.
\citet{ati2026} address this difficulty directly by coupling the
U-statistic to its Gaussian-chaos approximation and constructing
time-uniform boundaries for the resulting chaos process. This yields sharp
sequential calibration, but estimating the spectral quantities that
determine the Gaussian-chaos law involves eigendecomposition of empirical
Gram matrices, which can be computationally demanding under dense
monitoring. {
Most recently, \citet{LiuEtAl2026BootAVT} develop Boot-AVT and its
multi-kernel extension, using a quadratic multiplier bootstrap to
approximate the limiting Gaussian-chaos process without explicit
spectral estimation. 
}

One way to avoid this obstruction is to introduce a statistic or stochastic
structure with a more tractable sequential law. Cross U-statistics use
sample splitting and self-normalization to recover Gaussian calibration
\citep{KimRamdas2024}; this idea has led to permutation-free versions of
MMD, HSIC, and distance covariance
\citep{ShekharKimRamdas2022,ShekharKimRamdas2023}.
Martingale and betting methods provide a different route to sequential
nonparametric inference, including two-sample testing
\citep{ShekharRamdas2024}, kernelized independence testing
\citep{PodkopaevEtAl2023}, interactive rank and conditional-independence
testing \citep{DuanRamdasWasserman2022,GrunwaldHenziLardy2024}, and
sequential kernelized goodness-of-fit testing
\citep{MartinezTaboadaRamdas2025}.

More recently, \citet{ChatterjeeRamdas2025mMMD} introduce a
row-weighted MMD statistic with an asymptotically Gaussian null
distribution, yielding fixed-sample tests without spectral estimation.
Their multi-kernel extension also yields a Mahalanobis-type test based
on a chi-squared limit. These fixed-sample results do not by themselves
provide time-uniform error control under repeated monitoring. Building on this row-weighted construction, we develop asymptotic
time-uniform inference for general degenerate degree-two kernels. At each
time point, the pairwise terms involving the new observation are averaged
over the preceding sample and then accumulated across rows. This
reweighting preserves the target parameter, while degeneracy makes the
centered row averages a martingale-difference sequence. The resulting
martingale evolves on the harmonic variance scale
$T_n=H_{n-1}\asymp\log n$. The row normalization changes the variance from $O(n^{-2})$ for the
ordinary degenerate U-statistic to $O(\log n/n^2)$. Despite this logarithmic
inflation, the resulting confidence sequences remain nearly $n^{-1}$:
their radii are of order
$\sqrt{\log n\,\log\log\log n}/n$ for the stitched LIL boundary and
$\sqrt{\log n\,\log\log n}/n$ for the GM boundary. At the same time, the only nuisance quantity to be estimated is a scalar
variance, rather than the spectrum of the kernel operator.

We establish a strong Gaussian approximation on the harmonic clock and control the approximation and scale-estimation errors uniformly over an unbounded monitoring horizon. Combining these results with Gaussian-mixture (GM) and stitched law-of-the-iterated-logarithm (LIL) boundaries yields feasible confidence sequences with asymptotic time-uniform coverage and corresponding one-sided sequential tests. The
scalar intervals also form asymptotic confidence sequences in the sense of
\citet{WaudbySmithArbourSinhaKennedyRamdas2024}. The same framework extends
to multi-dimensional vector-valued kernels: a strong Gaussian approximation
for the vector martingale and consistent covariance estimation lead to
Mahalanobis confidence ellipsoids calibrated by radial GM and LIL-type
Brownian boundaries
\citep{de_la_Pe_a_2007,WhitehouseWuRamdas2025}.

Numerical experiments examine the finite-sample behavior of the resulting
scalar and vector procedures. For MMD and HSIC, we compare the one-sided
sequential tests with SAGE \citep{ati2026} as a spectral-calibration
baseline, focusing on sequential type-I error control and power. The vector
experiments assess the same metrics for three-sample homogeneity testing over the monitoring windows considered.

The rest of the paper is organized as follows.
Section~\ref{sec:methodology} presents the proposed method and
representative examples. Section~\ref{sec:scalar-results} establishes
the main theoretical results, and Section~\ref{sec:vector-extension}
develops an extension to multi-dimensional kernels.
Section~\ref{app:online-implementation} describes online computation.
Section~\ref{sec:simulations} presents simulations, and
Section~\ref{sec:discussion} concludes.
Proofs and their supporting lemmas are collected in the appendix.

\section{Methodology}\label{sec:methodology}

Let $X_1,X_2,\ldots$ be i.i.d. random elements on
$\left(\calX,\calA\right)$ with distribution $P$.
Let $h:\calX^2\to\R$ be a symmetric measurable kernel, and let
$\theta=\E h(X_1,X_2)$ be the parameter of interest.
Write $h_0(x,y)=h(x,y)-\theta$.
We assume that the centered kernel is first-order degenerate
(also called canonical), $ \E\left\{h_0(X_1,x)\right\}=0 $ for $P\text{-a.e. }x.$ The ordinary degree-two U-statistic is
\(U_n=\binom n2^{-1}\sum_{1\le i<j\le n}h(X_i,X_j).\)
Under square integrability and degeneracy,
$n(U_n-\theta)$ has a quadratic Gaussian-chaos limit.
We instead change the weighting of the rows while retaining the target
$\theta$. For $n\ge2$, define
\begin{align*}
    \widetilde U_n
  &=\frac1{n-1}\sum_{i=2}^n \frac1{i-1}\sum_{j<i}h(X_i,X_j) \\
  &=\frac1{n-1}\left[\sum_{i=2}^n \frac1{i-1}
      \sum_{j<i}\left\{h(X_i,X_j)-\theta\right\}\right]+\theta
    :=\frac{M_n}{n-1}+\theta.
\end{align*}
Degeneracy implies that $M_n$ is a martingale with respect to
$\mathcal F_n=\sigma(X_1,\ldots,X_n)$. If
$\tau^2:=\E h_0^2(X_1,X_2)<\infty$, then the centered row averages
$Y_i=(i-1)^{-1}\sum_{j<i}h_0(X_i,X_j)$ satisfy
$\E Y_i^2=\tau^2/(i-1)$. Write
$T_n=H_{n-1}=\sum_{r=1}^{n-1}r^{-1}$ and set $M_1=0$.
The resulting orthogonality gives $ \Var(\widetilde U_n)={\tau^2T_n}/{(n-1)^2},$ while $\Var(U_n)={2\tau^2}/{n(n-1)}$. Thus row normalization incurs only a
logarithmic variance inflation.

The martingale construction above applies whenever the centered kernel is
degenerate. This condition arises naturally under the null hypothesis in a
number of kernel- and distance-based testing problems. We record three
representative examples, two of which will also be used in the numerical
experiments. 

\noindent\textbf{Maximum mean discrepancy (MMD).}
Let $Z_i=\left(X_i,Y_i\right)$, where $X_i\sim P$, $Y_i\sim Q$, and the pairs
$Z_1,Z_2,\ldots$ are i.i.d. Let $k:\mathcal X^2\to\mathbb R$ be a
symmetric positive definite kernel. The squared maximum mean discrepancy is
estimated by the kernel
\(
  h(z,z')
  =
  k(x,x')+k(y,y')-k(x,y')-k(x',y).
\)
Its target parameter is
$\theta=\E h(Z_1,Z_2)=\operatorname{MMD}^2_k(P,Q)$.
Under the null hypothesis $H_0:P=Q$, $\theta=0$ and
$\E\left\{h(z,Z_1)\right\}=0$ for $P_Z$-almost every $z$, so the kernel is
degenerate. For bounded kernels, such as the Gaussian kernel $k(x,y)=\exp\{-\|x-y\|^2/(2\sigma^2)\}$, the
fourth-moment condition imposed in Section~\ref{sec:scalar-results} is immediate.

\medskip

\noindent\textbf{Hilbert--Schmidt independence criterion (HSIC).}
Let $Z_i=\left(X_i,Y_i\right)\sim P_{XY}$ be i.i.d. observations on
$\mathcal X\times\mathcal Y$. Let $k$ and $\ell$ be symmetric positive
definite kernels on $\mathcal X$ and $\mathcal Y$, respectively. Define the
population-centered kernels $\widetilde k(x,x')
  = k(x,x')-\E k(x,X_1)-\E k(X_1,x')+\E k(X_1,X_2)$ and $\widetilde \ell(y,y')
  = \ell(y,y')-\E \ell(y,Y_1)-\E \ell(Y_1,y')+\E \ell(Y_1,Y_2).$ Then the population-centered HSIC kernel is
$h(z,z')=\widetilde k(x,x')\widetilde \ell(y,y')$.
Its target parameter is
$\theta=\E h(Z_1,Z_2)=\operatorname{HSIC}(P_{XY},P_XP_Y)$.
Under the null hypothesis $H_0:X\perp Y$, $\theta=0$, and
\(\E\left\{h(z,Z_1)\right\} = \E\left\{\widetilde k(x,X_1)\widetilde \ell(y,Y_1)\right\} = \E\left\{\widetilde k(x,X_1)\right\} \E\left\{\widetilde \ell(y,Y_1)\right\} = 0.\)
Thus $h_0=h$ is degenerate under $H_0$. If $k$ and $\ell$ are bounded,
then $h$ is bounded. This construction uses population centering; replacing the population
quantities by empirical estimates requires additional analysis beyond
the fixed-kernel setting considered here.

\medskip

\noindent\textbf{Quadratic empirical-process functionals.}
A broader class is obtained from quadratic empirical-process
discrepancies. Let $\left\{\phi_t:t\in\mathcal T\right\}$ be a measurable
class of real-valued functions, and let $\nu$ be a finite measure on
$\mathcal T$. Suppose the null hypothesis specifies the moment function $m_0(t)=\E_0\left\{\phi_t(X)\right\},$ $t\in\mathcal T.$
Assume joint measurability in $(t,x)$ and
$\E\int |\phi_t(X)-m_0(t)|^2\,d\nu(t)<\infty$, which permit the expectation
and integration steps below. A quadratic empirical-process discrepancy is
\(
  \theta
  =
  \int_{\mathcal T}
  \left\{\E\phi_t(X)-m_0(t)\right\}^2\,d\nu(t).
\)
It is estimated by the degree-two kernel
\(h(x,y) = \int_{\mathcal T} \left\{\phi_t(x)-m_0(t)\right\} \left\{\phi_t(y)-m_0(t)\right\}\,d\nu(t).\)
Indeed, $\E h(X_1,X_2)
  =
  \int_{\mathcal T}
  \left\{\E\phi_t(X)-m_0(t)\right\}^2\,d\nu(t)
  =
  \theta.$ Under the null hypothesis
$H_0:\E\phi_t(X)=m_0(t)$ for all $t\in\mathcal T$, we have $\theta=0$, and
$\E\left\{h(x,X_1)\right\}=0.$
Hence $h_0=h$ is degenerate under $H_0$. Classical Cram\'er--von Mises
statistics and related weighted quadratic goodness-of-fit discrepancies
are included in this class.

The preceding examples illustrate how the degenerate condition arises under the null hypothesis. To formulate the corresponding sequential guarantees, we recall two related definitions of asymptotic confidence sequences
discussed by \citet{WaudbySmithArbourSinhaKennedyRamdas2024}.
The first concerns asymptotic time-uniform coverage, while the second
characterizes asymptotic equivalence to an exact confidence sequence.

\begin{definition}[Asymptotic time-uniform coverage]
\label{def:delayed-coverage}
For each deterministic integer $m\ge2$, let $\{C_{m,n}\}_{n\ge m}$ be a
sequence of random sets. The family is said to have asymptotic
time-uniform $(1-\alpha)$-coverage for $\theta$ if
\[
\liminf_{m\to\infty}
\Pr\!\left(
\theta\in C_{m,n}
\text{ for every } n\ge m
\right)
\ge 1-\alpha.
\]
\end{definition}

Definition~\ref{def:delayed-coverage} requires simultaneous coverage over
the entire tail $n\ge m$, with the coverage guarantee obtained as the
deterministic monitoring start $m$ tends to infinity.

\begin{definition}[Asymptotic confidence sequences]
\label{def:comparison}
For a fixed starting point $m$, intervals
$C_{m,n}=[\widehat\theta_n\pm w_{m,n}]$, $n\ge m$, form an asymptotic
confidence sequence if there exist intervals
$C^*_{m,n}=[\widehat\theta_n\pm w^*_{m,n}]$ satisfying
\[
 \Pp\{\theta\in C^*_{m,n}\text{ for every }n\ge m\}\ge1-\alpha,
 \qquad
 {w^*_{m,n}}/{w_{m,n}}\longrightarrow1
 \quad\text{almost surely as }n\to\infty.
\]
\end{definition}

Definition~\ref{def:comparison} instead fixes $m$ and requires the
candidate intervals to be asymptotically equivalent in width to an exact
confidence sequence. These two properties are distinct: asymptotic equivalence to an exact
confidence sequence at a fixed starting point does not by itself imply
the time-uniform coverage property in
Definition~\ref{def:delayed-coverage}.

We next record the Gaussian boundary results used in the
subsequent analysis. In the scalar case, these are expressed through time-uniform
boundaries for standard Brownian motion.

\begin{definition}[Brownian motion boundary]
\label{def:scalar-boundary}
For $t_0>0$, a nondecreasing function
$u_{\alpha,t_0}:\left[t_0,\infty\right)\to\left(0,\infty\right)$ is a
two-sided $(1-\alpha)$ Brownian motion boundary if, for standard Brownian
motion $W$,
\(\Pp\left\{ \exists t\ge t_0:|W(t)|\ge u_{\alpha,t_0}(t) \right\} \le\alpha .\)
\end{definition}

Two standard choices are used throughout.

\noindent\textbf{Robbins Gaussian-mixture boundary.}
Let $g(a)=2\{1-\Phi(a)+a\phi(a)\}$ for $a\ge0$, where $\Phi$ and $\phi$ are
the standard normal distribution and density functions. The scalar
Gaussian-mixture boundary may be written as
\begin{equation}
  u_{\alpha,t_0}^{\rm GM}(t)
  =
  \sqrt{t\left[\left\{g^{-1}(\alpha)\right\}^2
  +\log\left(t/t_0\right)\right]},
  \qquad t\ge t_0 .
  \label{eq:gm-boundary}
\end{equation}

\noindent\textbf{Stitched LIL boundary.}
For $\eta>1$ and $s>1$,
\begin{equation}
  u_{\alpha,t_0}^{\rm LIL}(t)
  =
  \frac{\eta^{1/4}+\eta^{-1/4}}{\sqrt2}
    \sqrt{t\left[
      s\log\log\left(\max\left\{\frac{\eta t}{t_0},e\right\}\right)
      +\log\left\{\frac{2\zeta(s)}
      {\alpha\left(\log\eta\right)^s}\right\}
    \right]} .
  \label{eq:lil-boundary}
\end{equation}
Here $\zeta$ is the Riemann zeta function. The one-sided stitched formula
of \citet{HowardRamdasMcAuliffeSekhon2021} has $\zeta(s)$ in place of
$2\zeta(s)$. We allocate $\alpha/2$ to each sign in
\eqref{eq:lil-boundary}.

Lemma~\ref{lem:scalar-brownian-calibration}
establishes the two-sided validity of both displayed boundaries. For fixed $t_0$, the LIL boundary has order
$O(\sqrt{t\log\log t})$, whereas the Gaussian-mixture boundary has order
$O(\sqrt{t\log t})$.
\section{Main Results}\label{sec:scalar-results}

Throughout this section, we retain the setup and degeneracy condition of
Section~\ref{sec:methodology} and impose the following additional condition:
\begin{itemize}
    \item[(A1)] The kernel satisfies $\E\left|h_0(X_1,X_2)\right|^4<\infty,$ and
    $\tau^2:=\E h_0^2(X_1,X_2)>0$.
\end{itemize}

We estimate the unknown scale \(\tau^2\) by the rowwise plug-in estimator
\begin{equation*}
  \widehat\tau_n^2
  =
  \frac1{T_n}\sum_{i=2}^n\left\{\frac1{i-1}\sum_{j<i}h(X_i,X_j)-\widetilde U_n\right\}^2.
\end{equation*}

{Recursive implementations of the scalar and vector procedures, together
with their computational costs, are given in
Section~\ref{app:online-implementation}.}
The following theorem supplies both its
consistency and the variance-clock control needed for the strong approximation.

\begin{theorem}\label{thm:scalar-clock}
Under (A1), $(Y_i,\calF_i)$ is a martingale-difference sequence with the following properties.
\begin{enumerate}[label=(\alph*)]
  \item The predictable quadratic variation satisfies
  \({\bracket{M}_n}/{T_n} :=\sum_{i=2}^n\E\left(Y_i^2\mid\calF_{i-1}\right) /{T_n}\asto \tau^2.\)
    Moreover, ${\bracket{M}_n}-\tau^2{T_n}= o_{a.s.}\left(T_n/\log T_n\right)$.
  \item The optional quadratic variation satisfies
  \({\qv{M}_n}/{T_n} :=\sum_{i=2}^nY_i^2/{T_n} \asto \tau^2 .\)
  \item The plug-in estimator is strongly consistent, that is $\widehat\tau_n^2\asto\tau^2$.
\end{enumerate}
\end{theorem}

Since $T_n\asymp\log n$, Theorem~\ref{thm:scalar-clock} identifies the
harmonic clock as the natural variance scale of $M_n$. The next theorem
establishes a strong Gaussian approximation on this time scale.

\begin{theorem}
\label{thm:scalar-asip}
Under (A1), on an enriched probability space there exists a standard Brownian motion $W$ such that
$M_n=\tau W(T_n)+ o_{a.s.}\left(\sqrt{T_n}\right).$
Equivalently,
\begin{equation*}
  \widetilde U_{n}-\theta
  =\frac{\tau W(T_n)}{n-1}+o_{a.s.}\left(\frac{\sqrt{T_n}}{n}\right).
\end{equation*}
\end{theorem}

Combining the Brownian approximation with either boundary in
Section~\ref{sec:methodology} gives the following sequence of intervals:
\begin{equation}
  C_{m,n}^{\rm MG}(\alpha)
  =
  \left[
    \widetilde U_n
    \pm
    \widehat\tau_n
    \frac{u_{\alpha,T_m}(T_n)}{n-1}
  \right],
  \qquad n\ge m .
  \label{eq:martingale-CS}
\end{equation}
The following theorem establishes asymptotic time-uniform coverage of \eqref{eq:martingale-CS} in the
sense of Definition~\ref{def:delayed-coverage}.

\begin{theorem}
\label{thm:scalar-asympcs}
Under (A1), let $u_{\alpha,t_0}$ be either \eqref{eq:gm-boundary} or \eqref{eq:lil-boundary}. Then
\begin{equation*}
  \liminf_{m\to\infty}
  \Pp\left(\forall n\ge m:\theta\in C_{m,n}^{\rm MG}(\alpha)\right)
  \ge1-\alpha .
\end{equation*}
\end{theorem}

Since $T_n\asymp\log n$ and $\widehat\tau_n\asto\tau$, the interval-width
orders follow directly from the growth rates of the two Brownian boundaries.
For fixed $m$ and the boundary tuning parameters,  the radii of $C_{m,n}^{\rm MG,LIL}$ and $C_{m,n}^{\rm MG,GM}$ are of order $
  O_{a.s.}\left({\sqrt{\log  n\,\log\log\log n}}/{n}
  \right)$ and $
  O_{a.s.}\left({\sqrt{\log n\log\log n}}/{n}
  \right)$ respectively.

\begin{remark}
\label{rem:scalar-asympcs}
For every fixed $m\ge2$, the sequence
$\{C^{\mathrm{MG}}_{m,n}(\alpha):n\ge m\}$ also satisfies Definition~\ref{def:comparison}. More precisely, under
(A1) and for either boundary, there exists a nonasymptotic
$(1-\alpha)$-confidence sequence
$\{C^*_{m,n}(\alpha):n\ge m\}$, centered at $\widetilde U_n$, such that
\(
   {\operatorname{rad}\{C^*_{m,n}(\alpha)\}}/
       {\operatorname{rad}\{C^{\mathrm{MG}}_{m,n}(\alpha)\}}
  \asto 1\), as \(n\to\infty.
\)
This property complements the time-uniform coverage guarantee
of Theorem~\ref{thm:scalar-asympcs}.
\end{remark}

\begin{remark}
\label{rem:one-sided-testing}
In many applications of degenerate U-statistics, degeneracy arises under a
null hypothesis and the relevant alternative is one-sided. Examples include
MMD and HSIC, for which the population discrepancy is nonnegative and the
kernel is degenerate under the null hypothesis with $\theta_0=0$.

For $\eta>1$, $s>1$, and $\alpha\in(0,1)$, the upper stitched LIL boundary is
\[
u_{\alpha,t_0}^{+,\mathrm{LIL}}(t)
=
\frac{\eta^{1/4}+\eta^{-1/4}}{\sqrt{2}}
\sqrt{
t\left[
s\log\log\left(
\max\left\{\frac{\eta t}{t_0},e\right\}
\right)
+
\log\left\{
\frac{\zeta(s)}
{\alpha(\log\eta)^s}
\right\}
\right]
},
\qquad t\ge t_0,
\]
which differs from \eqref{eq:lil-boundary} only in that the full error
probability $\alpha$ is allocated to the upper tail.

For the Gaussian-mixture construction, let
$g_+(a)=1-\Phi(a)+a\phi(a)+\phi(a)^2/\Phi(a)$, $a\ge0$, and, for
$0<\alpha<g_+(0)=1/2+1/\pi$, define
$g_+^{-1}(\alpha)=\inf\{a\ge0:g_+(a)\le\alpha\}$.
The corresponding upper boundary is
\[
  u_{\alpha,t_0}^{+,\mathrm{GM}}(t)
  =
  \sqrt{t\left\{[g_+^{-1}(\alpha)]^2+\log(t/t_0)\right\}},
  \qquad t\ge t_0.
\]
Lemma~\ref{lem:scalar-one-sided-boundaries} establishes that both boundaries
satisfy
\(
\Pp\left\{
\exists t\ge t_0:
W(t)\ge u_{\alpha,t_0}^+(t)
\right\}
\le\alpha
\)
for their respective parameter ranges.
\end{remark}

\begin{corollary}
\label{cor:one-sided-testing}
Fix a distribution under $H_0:\theta=\theta_0$ for which the centered kernel
$h-\theta_0$ is degenerate and satisfies (A1). Let $u_{\alpha,t_0}^+$ be
either boundary in Remark~\ref{rem:one-sided-testing}, with parameters in
the stated range. Then
\[
\limsup_{m\to\infty}
\Pp_{H_0}\left\{
\exists n\ge m:
(n-1)(\widetilde U_n-\theta_0)
>
\widehat\tau_n
u_{\alpha,T_m}^+(T_n)
\right\}
\le\alpha.
\]
Thus, the event in the preceding display is the rejection region for testing
$H_0:\theta=\theta_0$ against $H_1:\theta>\theta_0$, with rejection occurring
at the first boundary crossing.
\end{corollary}

The following result imposes the fourth-moment condition in (A1)
under the alternative distribution, without assuming degeneracy.

\begin{proposition}[Fixed-alternative consistency]
\label{prop:one-sided-fixed-alternative}
Let $X_1,X_2,\ldots$ be i.i.d.\ with distribution $P$, and let $h$ be a fixed
symmetric measurable kernel with
$\theta=\E h(X_1,X_2)>\theta_0$.
Suppose that $\E|h_0(X_1,X_2)|^4<\infty$, where $h_0=h-\theta$.
Fix a deterministic integer $m\ge2$, and let $u^+_{\alpha,t_0}$ be
either one-sided boundary in Remark~\ref{rem:one-sided-testing},
with all parameters fixed within their stated ranges.
Then, almost surely,
\[
  (n-1)(\widetilde U_n-\theta_0)
  >\widehat\tau_n u^+_{\alpha,T_m}(T_n)
  \qquad
  \text{for all sufficiently large }n\ge m.
\]
In particular, the first rejection time is finite almost surely,
and the probability of rejection by time $N$ tends to one as
$N\to\infty$.
\end{proposition}

The proof is given in Appendix~\ref{app:one-sided-fixed-alternative}.

\section{Extension to multi-dimensional kernels}
\label{sec:vector-extension}

We now extend the martingaleized construction to multi-dimensional
kernels. Let $\bm h:\calX^2\to\R^d$ be a symmetric measurable
kernel, where $d$ is fixed, and define
$\bm\theta=\E\bm h(X_1,X_2)\in\R^d$ and
$\bm h_0(x,y)=\bm h(x,y)-\bm\theta$. We assume that the centered kernel is degenerate componentwise:
\(
  \E\left\{\bm h_0(X_1,x)\right\}=\bm 0\in\R^d\) for \(P\text{-a.e. }x.
\)
Define the corresponding statistic
\[
  \widetilde{\bm U}_n
  =
  \frac1{n-1}
  \sum_{i=2}^n
  \frac1{i-1}
  \sum_{j<i}\bm h(X_i,X_j),
\]
and, analogously to the scalar case, denote
$\bm Y_i=\sum_{j<i}\bm h_0(X_i,X_j)/(i-1)$ and
$\bm M_n=\sum_{i=2}^n\bm Y_i$. Then
$\bm M_n$ is an $\R^d$-valued martingale with respect to
$\calF_n=\sigma(X_1,\ldots,X_n)$.

Throughout this section, we impose the following additional condition:
\begin{itemize}
  \item[(A2)] The kernel satisfies
  $\E\norm{\bm h_0(X_1,X_2)}_2^4<\infty$, and
  $\bm\Omega:=\E\left\{
      \bm h_0(X_1,X_2)
      \bm h_0(X_1,X_2)\trans
    \right\}$ is positive definite.
\end{itemize}

We estimate $\bm\Omega$ by the rowwise plug-in estimator
\[
  \widehat{\bm\Omega}_n
  =
  \frac1{T_n}
  \sum_{i=2}^n
  \left\{
    \frac1{i-1}\sum_{j<i}\bm h(X_i,X_j)
    -
    \widetilde{\bm U}_n
  \right\}
  \left\{
    \frac1{i-1}\sum_{j<i}\bm h(X_i,X_j)
    -
    \widetilde{\bm U}_n
  \right\}\trans.
\]
Lemma~\ref{lem:vector-clock} establishes that
$\widehat{\bm\Omega}_n\asto\bm\Omega$ under (A2). We next give a strong
Gaussian approximation for the vector martingale on the harmonic clock.

\begin{proposition}
\label{prop:vector-asip}
Under (A2), on an extension of the underlying probability space,
there exist a standard \(d\)-dimensional Brownian motion \(\bm W\)
and a constant \(\kappa>0\) such that
\[
  \norm{
    \bm M_n-\bm\Omega^{1/2}\bm W(T_n)
  }_2
  =
  O_{a.s.}(T_n^{1/2-\kappa})
  =
  o_{a.s.}({T_n^{1/2}}).
\]
\end{proposition}

To translate this strong Gaussian approximation into joint time-uniform
inference, we use radial boundaries for the Euclidean norm of a
$d$-dimensional Brownian motion. These are the vector counterparts of the
scalar Brownian boundaries introduced in Section~\ref{sec:methodology}.

\begin{definition}[Radial Brownian motion boundary]\label{def:vector-boundary}
For $t_0>0$, a nondecreasing function
$b_{d,\alpha,t_0}:[t_0,\infty)\to(0,\infty)$ is called a radial
$(1-\alpha)$ Brownian motion boundary if, for a standard
$d$-dimensional Brownian motion $\bm W_d$,
\(
  \Pp\left\{
    \exists t\ge t_0:
    \norm{\bm W_d(t)}_2
    \ge
    b_{d,\alpha,t_0}(t)
  \right\}
  \le\alpha.
\)
\end{definition}

We use the following Gaussian-mixture and LIL-type radial boundaries.

\noindent\textbf{Gaussian-mixture radial boundary.}
An explicit radial boundary follows from a Gaussian mixture of Brownian exponential martingales; see \citet{de_la_Pe_a_2007}
and \citet[Fact~4.3]{WhitehouseWuRamdas2025}. For a standard
$d$-dimensional Brownian motion, the boundary is
\begin{equation}
    b^{\mathrm{GM}}_{d,\alpha,t_0}(t)
    =
    \left[
        (t+t_0)
        \left\{
            d\log\left(1+\frac{t}{t_0}\right)
            +2\log(1/\alpha)
        \right\}
    \right]^{1/2},
    \qquad t\ge t_0 .
    \label{eq:vector-gm-boundary}
\end{equation}

\noindent\textbf{LIL-type radial boundary.}
We use the vector self-normalized inequality of
\citet[Theorem~4.1]{WhitehouseWuRamdas2025}.
Fix \(a>1\), \(\beta>1\), \(\varepsilon\in(0,1)\), and \(q>1\).
For \(t\ge t_0\), define
\begin{equation}
\begin{aligned}
    b^{\mathrm{LIL}}_{d,\alpha,t_0}(t)
    =
    \frac{\sqrt{2at}}{1-\varepsilon}
    \Bigg[
        &q\log\left\{
            1+\log_a\left(\frac{t}{t_0}\right)
        \right\}
       +\log\frac{\zeta(q)\beta}{\alpha(1-\beta^{-1})}
        +d\log\left(\frac{3\beta}{\varepsilon}\right)
    \Bigg]^{1/2}.
\end{aligned}
\label{eq:vector-lil-boundary}
\end{equation}

For fixed \(t_0\), \(d\), and tuning parameters,
the Gaussian-mixture boundary has order
\(O(\sqrt{t\log t})\), whereas the LIL-type boundary has the
sharper order \(O(\sqrt{t\log\log t})\).
Both choices are valid radial Brownian boundaries.

Using the radial boundary $b_{d,\alpha,t_0}$, define, for
$n\ge m$, the Mahalanobis confidence ellipsoid
\begin{equation}
  C_{m,n}^{\rm ell}(\alpha)
  =
  \left\{
    \bm\vartheta\in\R^d:
    (n-1)^2
    \left(
      \widetilde{\bm U}_n-\bm\vartheta
    \right)\trans
    \widehat{\bm\Omega}_n^{-1}
    \left(
      \widetilde{\bm U}_n-\bm\vartheta
    \right)
    \le
    \left\{
      b_{d,\alpha,T_m}(T_n)
    \right\}^2
  \right\},
  \label{eq:vector-ellipsoid}
\end{equation}
with the convention that
$C_{m,n}^{\rm ell}(\alpha)=\R^d$ whenever
$\widehat{\bm\Omega}_n$ is singular.

The following theorem gives asymptotic time-uniform coverage of these
ellipsoids in the sense of Definition~\ref{def:delayed-coverage}.

\begin{theorem}
\label{thm:vector-coverage}
Suppose that (A2) holds, and choose either \eqref{eq:vector-gm-boundary} or \eqref{eq:vector-lil-boundary}. Then
\[
  \liminf_{m\to\infty}
  \Pp\left\{
    \forall n\ge m:
    \bm\theta\in C_{m,n}^{\rm ell}(\alpha)
  \right\}
  \ge1-\alpha.
\]
\end{theorem}

\begin{remark}
\label{rem:vector-asympcs}
Under (A2), for every fixed $m\ge2$ and either radial boundary, the
ellipsoids satisfy a multivariate analogue of Definition~\ref{def:comparison}. Specifically, in the sense of
\citet[Definition~B.6]{WaudbySmithArbourSinhaKennedyRamdas2024}, there exists a
nonasymptotic
\((1-\alpha)\)-confidence sequence
\(\{C_{m,n}^{*,\rm ell}(\alpha):n\ge m\}\), centered at
$\widetilde{\bm U}_n$, such that
\(   {
    \nu\left(
      C_{m,n}^{\rm ell}(\alpha)
      \mathbin{\triangle}
      C_{m,n}^{*,\rm ell}(\alpha)
    \right)
  }/{
    \nu\left(
      C_{m,n}^{*,\rm ell}(\alpha)
    \right)
  }
  \asto0\)  as \(
 n\to\infty.
\)
Here $\nu$ denotes Lebesgue measure on $\R^d$, and $\triangle$ denotes
symmetric difference. 
\end{remark}

\section{Online implementation and computational cost}
\label{app:online-implementation}

\paragraph{Scalar updates.}
For each new observation, compute the uncentered row average $ a_n=\sum_{j<n}h(X_n,X_j)/(n-1)$ where $ n\ge2.$ The statistic and its scale estimator can be written as $\widetilde U_n=\sum_{i=2}^n a_i/(n-1)$ and $\widehat\tau_n^2={S_n}/{T_n}=\sum_{i=2}^n(a_i-\widetilde U_n)^2/{T_n}.$

Algorithm~\ref{alg:online-mg} updates the mean and centered sum of
squares recursively. The identity
$S_n=\sum_{i=2}^n a_i^2-(n-1)\widetilde U_n^2$
verifies that these updates compute the original estimators exactly
in exact arithmetic. Maintaining the centered sum of squares also
avoids subtracting two potentially large raw second-moment terms.
All observations, including those received before the monitoring
start $m$, contribute to the updates.

\begin{algorithm}[ht]
\caption{Online scalar martingaleized inference}
\label{alg:online-mg}
\begin{algorithmic}[1]
\Require Kernel $h$, observations $X_1,X_2,\ldots$,
monitoring start $m\ge2$, and level $\alpha$.
\Require Selected two-sided boundary $u_{\alpha,t_0}$;
for one-sided testing, also select $u^+_{\alpha,t_0}$
and specify $\theta_0$.
\State Store $X_1$; initialize
$\widetilde U_1\gets0$, $S_1\gets0$, and $T_1\gets0$.
\For{$n=2,3,\ldots$}
  \State Receive $X_n$ and compute
  $a_n\gets(n-1)^{-1}\sum_{j<n}h(X_n,X_j)$.
  \State $\Delta\gets a_n-\widetilde U_{n-1}$.
  \State $\widetilde U_n\gets
  \widetilde U_{n-1}+\Delta/(n-1)$.
  \State $S_n\gets S_{n-1}
  +\{(n-2)/(n-1)\}\Delta^2$.
  \State $T_n\gets T_{n-1}+1/(n-1)$;
  $\widehat\tau_n^2\gets S_n/T_n$;
  $\widehat\tau_n\gets\sqrt{\widehat\tau_n^2}$.
  \If{$n=m$}
    \State Store $t_0\gets T_n$.
  \EndIf
  \If{$n\ge m$}
    \State Report
    $C^{\mathrm{MG}}_{m,n}(\alpha)
    =[\widetilde U_n
    \pm\widehat\tau_n u_{\alpha,t_0}(T_n)/(n-1)]$.
    \State For one-sided testing, reject at the first time satisfying
    \Statex \hspace{\algorithmicindent}
    $(n-1)(\widetilde U_n-\theta_0)>
    \widehat\tau_n u^+_{\alpha,t_0}(T_n)$.
  \EndIf
  \State Store $X_n$ for subsequent row evaluations.
\EndFor
\end{algorithmic}
\end{algorithm}

\paragraph{Vector updates.}
For a $d$-dimensional kernel, define
$\mathbf a_n=(n-1)^{-1}\sum_{j<n}\mathbf h(X_n,X_j)$
and
$\boldsymbol\Delta_n
=\mathbf a_n-\widetilde{\mathbf U}_{n-1}$.
Starting from $\widetilde{\mathbf U}_1=\mathbf0$ and
$\mathbf S_1=\mathbf0_{d\times d}$, we use the update
\begin{align*}
  \widetilde{\mathbf U}_n
  &=\widetilde{\mathbf U}_{n-1}
    +\frac{\boldsymbol\Delta_n}{n-1},\\
  \mathbf S_n
  &=\mathbf S_{n-1}
    +\frac{n-2}{n-1}
      \boldsymbol\Delta_n\boldsymbol\Delta_n^{\mathsf T}
  \ \ \text{and}\ \  %
  \widehat{\boldsymbol\Omega}_n
  =\frac{\mathbf S_n}{T_n}.
\end{align*}
These updates maintain
$\mathbf S_n=\sum_{i=2}^n
(\mathbf a_i-\widetilde{\mathbf U}_n)
(\mathbf a_i-\widetilde{\mathbf U}_n)^{\mathsf T}$,
and hence recover the covariance estimator in Section~\ref{sec:vector-extension}.
At monitored times, use the selected radial boundary to construct
the Mahalanobis ellipsoid. The ellipsoid is
$\mathbb R^d$ when $\widehat{\boldsymbol\Omega}_n$ is singular.
When it is positive definite, the Mahalanobis statistic can be
evaluated by a Cholesky solve.

\paragraph{Computational cost.}
At time $n$, computing the new row average requires $n-1$ kernel
evaluations. Given this average, updating the estimator, harmonic
clock, and scalar variance estimate takes $O(1)$ arithmetic
operations; the vector covariance update takes $O(d^2)$.
A fresh Cholesky factorization for the Mahalanobis statistic costs
$O(d^3)$ at each monitored time. Hence, for fixed $d$ and
constant-cost kernel evaluations, the total computational cost
through time $N$ is $O(N^2)$. The implementation retains past
observations and uses $O(d^2)$ additional storage, without storing
past kernel values.

{
The Boot-AVT method of \citet{LiuEtAl2026BootAVT} consistently
approximates the limiting Gaussian-chaos law under their assumptions.
For a direct implementation with a dense degree-two kernel, maintaining
$B$ quadratic multiplier trajectories adds $O(Bn)$ arithmetic per
observation and $O(BN^2)$ through time $N$. In the experiments in
Sections~\ref{sec:sim-vector} and~\ref{app:vector-coverage}, its
empirical coverage falls below the nominal level.}

\paragraph{Runtime comparison.}
We compare the per-update running times of the proposed MG procedure
and SAGE under the MMD and HSIC null models of
Section~\ref{sec:sim-scalar}. At
$n\in\{400,800,1200,1600,2000\}$,
Table~\ref{tab:scalar-update-runtime} reports the median time over ten
independent paths to update both the GM and LIL decisions.
Timings include kernel evaluations and, for SAGE, spectral computation,
but exclude data generation and initialization. All procedures were
timed on the same Apple M4. The $n=400$ spectral update uses a direct eigendecomposition of
a $200\times200$ matrix; the remaining checkpoints use a warm-started
subspace eigensolver.

\begin{table}[ht]
\centering
\caption{Median time in milliseconds for one monitoring
update under the null, based on ten independent paths.}
\label{tab:scalar-update-runtime}
\begin{tabular}{rcccc}
\toprule
& \multicolumn{2}{c}{MMD} & \multicolumn{2}{c}{HSIC} \\
\cmidrule(lr){2-3}\cmidrule(lr){4-5}
$n$ & MG & SAGE & MG & SAGE \\
\midrule
400  & 0.014 &  5.089 & 0.013 &  5.916 \\
800  & 0.020 &  3.532 & 0.018 &  6.289 \\
1200 & 0.027 &  6.959 & 0.024 & 12.325 \\
1600 & 0.034 & 12.646 & 0.029 & 22.761 \\
2000 & 0.043 & 19.197 & 0.035 & 33.734 \\
\bottomrule
\end{tabular}
\end{table}

\section{Simulation}
\label{sec:simulations}

We study the finite-sample behavior of the scalar and vector procedures
developed in Sections~\ref{sec:scalar-results} and
\ref{sec:vector-extension}. For the scalar case, we evaluate the one-sided
sequential tests in MMD and HSIC problems and compare them with the spectral
SAGE procedures of \citet{ati2026}. For the vector case, we examine simultaneous coverage and rejection behavior in a three-sample kernel testing problem. Throughout, the nominal level is $\alpha=0.05$.

\subsection{One-Sided Sequential Testing}
\label{sec:sim-scalar}

\paragraph{Setup and procedures.}
We consider the MMD and population-centered HSIC examples introduced in
Section~\ref{sec:methodology}. For MMD, let $Z_i=(X_i,Y_i)$, where
$X_i\sim \mathcal{N}(0,1)$ and $Y_i\sim N(\delta,1)$ are independent. For
$z=(x,y)$ and $z'=(x',y')$, use $k(x,x')=\exp\{-(x-x')^2/2\}$ and $h(z,z')
  =
  k(x,x')+k(y,y')-k(x,y')-k(x',y)$. The null corresponds to $\delta=0$, and the alternatives are
$\delta\in\{0.05,0.10,\ldots,0.45\}$.

For HSIC, we generate $(X_i,Y_i)$ from a centered bivariate normal
distribution with unit marginal variances and correlation $\rho$, and use
$h(z,z')=\widetilde k(x,x')\widetilde k(y,y')$, where the
population-centered Gaussian kernel is available analytically $\widetilde k(x,x')
  =
  e^{-(x-x')^2/2}
  -{(e^{-x^2/4}+e^{-(x')^2/4})}/{\sqrt2}
  +1/{\sqrt3}.$
The independence null corresponds to $\rho=0$, and the alternatives are
$\rho\in\{0.05,0.10,\ldots,0.50\}$. Both kernels are bounded and satisfy
(A1).

We implement the one-sided MG-GM and MG-LIL tests from
Section~\ref{sec:scalar-results} using the harmonic clock $T_n$ and the
rowwise plug-in variance estimator. MG-GM uses the one-sided
Gaussian-mixture boundary, while MG-LIL uses the one-sided stitched LIL
boundary with $\eta=2$ and $s=1.4$.

As spectral benchmarks, we use the one-sided SAGE-GM and SAGE-LIL
procedures of \citet{ati2026}, based on the ordinary U-statistic.
Following their simulation implementation, we retain
$L_n=\lfloor n^{1/4}\rfloor$ empirical eigenvalues and update the
spectral allocation weights proportionally to
$(\widehat\lambda_{\ell,n})_+$ at each reporting time $n$.
For SAGE-LIL, we use the same stitching parameters
as for MG-LIL.

\begin{figure*}[t]
  \centering
  \includegraphics[width=0.8\textwidth]
  {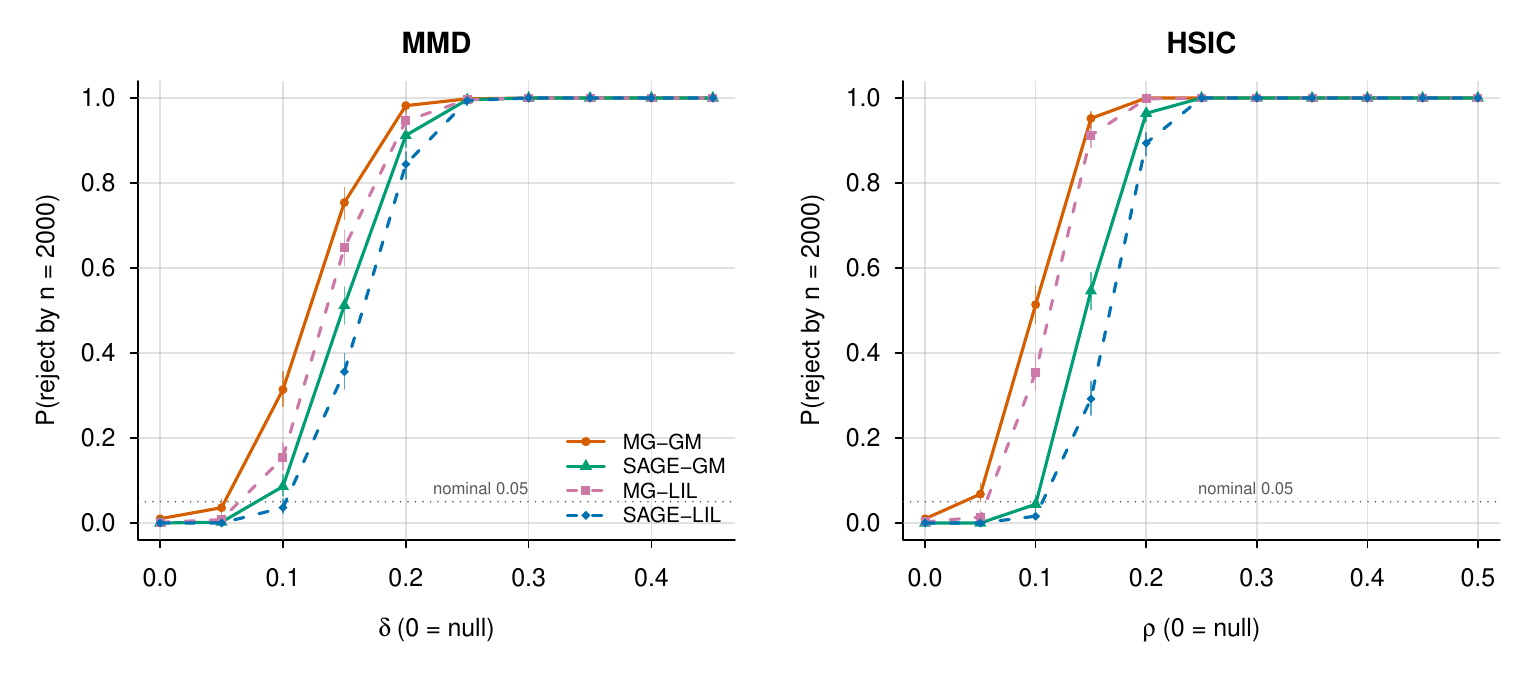}
  \caption{Cumulative rejection probability through $n=2000$ for the
  one-sided scalar tests. The left panel varies the MMD mean shift $\delta$,
  and the right panel varies the Gaussian correlation $\rho$ for HSIC.}
  \label{fig:scalar-one-sided-power}
\end{figure*}

We begin monitoring at $m=400$ over the grid
$\mathcal G=\{400,401,\ldots,2000\}$. For each example, empirical type-I error is evaluated using 2000 Monte
Carlo replications, and power at each nonzero effect size is evaluated
using 500 replications. Table~\ref{tab:scalar-sage} reports cumulative type-I error through
$n=2000$, while Figure~\ref{fig:scalar-one-sided-power} presents the
corresponding power curves. To compare detection speed at effect size
$0.15$, let
\(\nu=\inf\{n\in\mathcal G:
  \text{the rejection criterion is satisfied at time }n\},\)
with $\inf\varnothing=\infty$. Table~\ref{tab:scalar-first-rejection}
reports the empirical rejection probability by $N$ and the mean capped
rejection time $\min\{\nu,N+1\}$; a path without rejection by $N$ is
assigned $N+1=2001$ for this summary.

\begin{table}[!htbp]
  \centering
  \caption{Cumulative type-I error under monitoring over the grid
  $\mathcal G=\{400,401,\ldots,2000\}$.}
  \label{tab:scalar-sage}
  \small
  \begin{tabular}{lcccc}
    \toprule
    Example & MG-GM & SAGE-GM & MG-LIL & SAGE-LIL \\
    \midrule
    RBF-MMD & 0.0100 & 0.0000 & 0.0015 & 0.0000 \\
    HSIC    & 0.0100 & 0.0000 & 0.0025 & 0.0000 \\
    \bottomrule
  \end{tabular}
\end{table}

\paragraph{Results.}
Table~\ref{tab:scalar-sage} shows that the martingale procedures are less
conservative than SAGE while retaining empirical type-I error below the
nominal level; neither SAGE procedure rejects under the null. Moreover, the martingale procedures have substantially higher empirical power against moderate
alternatives, especially for HSIC: Table~\ref{tab:scalar-first-rejection}
gives the powers at effect size $0.15$. By effect size $0.25$, all four
procedures have power close to one, as shown in Figure~\ref{fig:scalar-one-sided-power}.

The mean capped first rejection times in
Table~\ref{tab:scalar-first-rejection} reinforce the power comparison.
At effect size $0.15$, the mean capped times for MG-GM and MG-LIL are
approximately $367$ and $306$ observations lower than those for
SAGE-GM and SAGE-LIL, respectively, for MMD. The corresponding
differences for HSIC are approximately $794$ and $779$ observations. Together with the higher endpoint power, these shorter detection times
provide additional evidence that the martingale procedures detect the
alternatives more rapidly in the settings considered.

The proposed procedures also offer substantial computational advantages
over SAGE; see Section~\ref{app:online-implementation} for timing comparisons.

\begin{table}[!htbp]
  \centering
  \caption{Rejection probability and mean capped first-rejection time
  at effect size $0.15$.}
  \label{tab:scalar-first-rejection}
  \small
  \begin{tabular}{llcccc}
    \toprule
    Example & Measure & MG-GM & SAGE-GM & MG-LIL & SAGE-LIL \\
    \midrule
    RBF-MMD, $\delta=0.15$
      & Rejection by $N$ & 0.754 & 0.512 & 0.648 & 0.356 \\
      & Mean capped time & 1239.8 & 1606.9 & 1460.6 & 1766.1 \\
    HSIC, $\rho=0.15$
      & Rejection by $N$ & 0.952 & 0.546 & 0.912 & 0.292 \\
      & Mean capped time & 881.9 & 1675.6 & 1068.8 & 1847.9 \\
    \bottomrule
  \end{tabular}
\end{table}

\subsection{Three-Sample Kernel Testing with Vector Confidence Ellipsoids}
\label{sec:sim-vector}

\paragraph{Setup and procedures.}
We consider the three-sample homogeneity null $H_0:P_1=P_2=P_3$.
At each time $i$, we observe
$\bm Z_i=(\bm X_i^{(1)},\bm X_i^{(2)},\bm X_i^{(3)})$,
where $\bm X_i^{(a)}\in\mathbb R^2$ has distribution $P_a$.
Observations are independent across groups and times. We use the fixed
Gaussian kernel
$k(\bm x,\bm y)=\exp\{-\|\bm x-\bm y\|_2^2/4\}$.
For $\bm z=(\bm x^{(1)},\bm x^{(2)},\bm x^{(3)})$,
$\bm z'=(\bm y^{(1)},\bm y^{(2)},\bm y^{(3)})$, and
$(a,b)\in\{(1,2),(1,3),(2,3)\}$, define $h_{ab}(\bm z,\bm z')
  =k(\bm x^{(a)},\bm y^{(a)})
   +k(\bm x^{(b)},\bm y^{(b)})
   -k(\bm x^{(a)},\bm y^{(b)})
   -k(\bm x^{(b)},\bm y^{(a)}).$
The vector kernel $\bm h=(h_{12},h_{13},h_{23})^\top$ has target $  \bm\theta
  =\E\{\bm h(\bm Z_1,\bm Z_2)\}
  =\bigl(
      \mathrm{MMD}_k^2(P_1,P_2),
      \mathrm{MMD}_k^2(P_1,P_3),
      \mathrm{MMD}_k^2(P_2,P_3)
    \bigr)^\top.$

Under $H_0$, $\bm\theta=\bm0$ and $\bm h$ is componentwise
degenerate. Under the null, we set $P_1=P_2=P_3=\mathcal N(\bm0,\mathbf I_2)$ and corresponding
$\bm\Omega=(3\mathbf I_3+\bm1_3\bm1_3^\top)/20$ is positive
definite. We set $P_1=P_2=\mathcal N(\bm0,\mathbf I_2)$
and $P_3=\mathcal N(\bm\mu_\delta,\mathbf I_2)$, where
$\bm\mu_\delta=(\delta,0)^\top$ and
$\delta\in\{0.1,0.2,0.3,0.4,0.5,0.6\}$ under the alternatives.
We use 2000 Monte Carlo replications under the null and 500 at
each nonzero effect size. Null coverage is evaluated for
$m\in\{20,30,40\}$, while the power curves use $m=30$; all
monitoring ends at $N=2000$.

We evaluate the ellipsoids in \eqref{eq:vector-ellipsoid} using
\eqref{eq:vector-gm-boundary} (VEC-GM) and
\eqref{eq:vector-lil-boundary} (VEC-WWR-LIL), with $d=3$,
$\alpha=0.05$, and $t_0=T_m$. For VEC-WWR-LIL, we set
$a=2$, $\beta=2$, $\varepsilon=0.5$, and $q=1.4$.
Here we report simultaneous coverage over $n=m,\ldots,N$ under
$H_0$ and the power across different $\delta$.

\begin{figure}[!htbp]
  \centering
  \begin{minipage}[b]{0.46\linewidth}
    \centering
    \begin{minipage}[c][3cm][c]{\linewidth}
      \centering
      \small
      \renewcommand{\arraystretch}{1.45}
      \setlength{\tabcolsep}{2pt}
      \begin{tabular}{rccc}
        \toprule
        $m$ & VEC-GM & VEC-WWR-LIL & Boot-AVT \\
        \midrule
        20 & 0.9620 & 1.0000 & {0.8750} \\
        30 & 0.9780 & 1.0000 & {0.8980} \\
        40 & 0.9810 & 1.0000 & {0.9125} \\
        \bottomrule
      \end{tabular}
    \end{minipage}
    \captionof{table}{Simultaneous coverage under the
    three-sample homogeneity null.}
    \label{tab:vector-coverage}
  \end{minipage}\hspace{0.06\linewidth}%
  \begin{minipage}[b]{0.46\linewidth}
    \centering
    \begin{minipage}[c][3cm][c]{\linewidth}
      \centering
      \includegraphics[width=\linewidth,height=4cm,keepaspectratio]
        {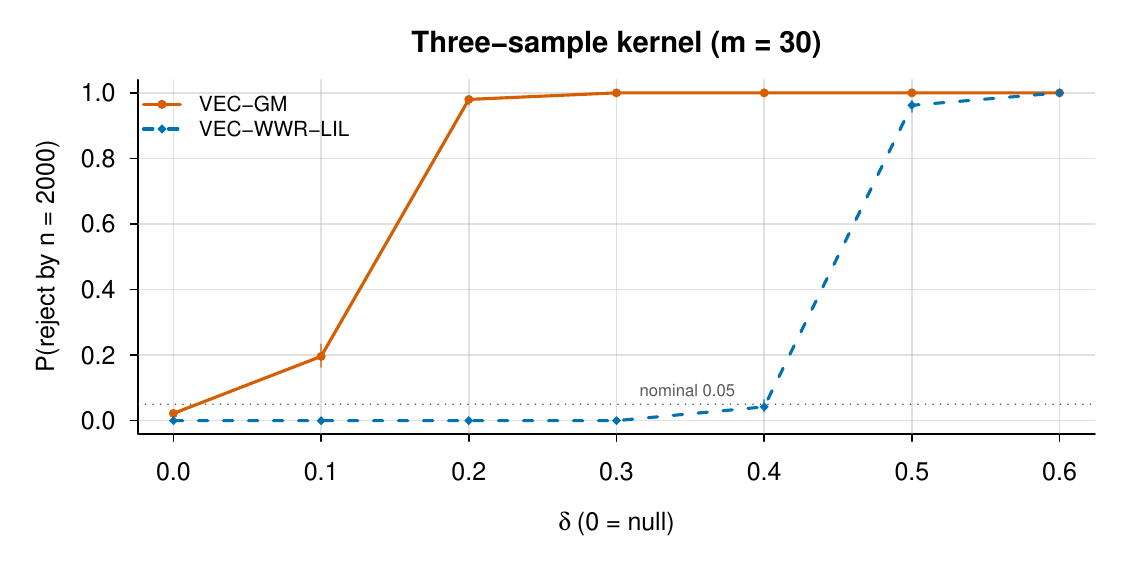}
    \end{minipage}
    \captionof{figure}{Rejection probability by $N=2000$
    against mean shift $\delta$, with $m=30$.}
    \label{fig:vector-three-sample-power}
  \end{minipage}
\end{figure}

\paragraph{Results.}
Table~\ref{tab:vector-coverage} shows empirical simultaneous
coverage between $0.9620$ and $0.9810$ for VEC-GM. VEC-WWR-LIL
covers all 2000 null paths at each monitoring start.
Figure~\ref{fig:vector-three-sample-power} shows a substantial
difference in detection power: VEC-GM rejects on nearly all paths
at $\delta=0.2$, whereas VEC-WWR-LIL rarely rejects through
$\delta=0.4$. The two procedures use the same estimator and
estimated covariance, so this difference reflects their radial
boundaries. At $n=2000$, the GM and LIL boundary values are
approximately $10.65$ and $42.76$, respectively. VEC-GM is thus
the more useful choice for detection within the monitoring window
considered.
{
Boot-AVT-Multi \citep{LiuEtAl2026BootAVT}, applied to the scaled
pairwise kernels $h_{ab}/\sqrt{0.4}$, falls short of the nominal
coverage at every monitoring start in Table~\ref{tab:vector-coverage}.
}

\subsection{Vector Product-Kernel Coverage}
\label{app:vector-coverage}

We also evaluate the vector ellipsoids for the kernel
$\bm h(x,y)=\mathbf A(x\odot y)$, where
$\mathbf A=\operatorname{diag}(1,2,1/2)$ and
$\odot$ denotes elementwise multiplication. The observations have
mean zero and covariance matrix $\mathbf\Sigma$ with unit diagonal
entries and off-diagonal entries equal to $1/2$. Consequently,
$\bm\theta=\bm0$ and
$\bm\Omega=\mathbf A(\mathbf\Sigma\odot\mathbf\Sigma)\mathbf A$.
We consider both $\mathcal N(\bm0,\mathbf\Sigma)$ observations and
a standardized multivariate $t_5$ distribution with covariance
$\mathbf\Sigma$. For each distribution, 2000 Monte Carlo paths are
monitored from $m\in\{20,30,40\}$ through $N=2000$, using the
VEC-GM and VEC-WWR-LIL procedures specified in
Section~\ref{sec:sim-vector}.

\begin{table}[!htbp]
  \centering
  \caption{Simultaneous coverage for the vector product kernel under
  Gaussian and standardized $t_5$ observations.}
  \label{tab:vector-product-coverage}
  \small
  \setlength{\tabcolsep}{3pt}
  \begin{tabular}{rcccccc}
    \toprule
    & \multicolumn{3}{c}{Gaussian}
    & \multicolumn{3}{c}{$t_5$} \\
    \cmidrule(lr){2-4}\cmidrule(lr){5-7}
    $m$
    & VEC-GM & VEC-WWR-LIL & Boot-AVT
    & VEC-GM & VEC-WWR-LIL & Boot-AVT \\
    \midrule
    20 & 0.9640 & 1.0000 & {0.8705} & 0.9750 & 1.0000 & {0.8660} \\
    30 & 0.9795 & 1.0000 & {0.8895} & 0.9855 & 1.0000 & {0.8855} \\
    40 & 0.9875 & 1.0000 & {0.9065} & 0.9885 & 1.0000 & {0.9005} \\
    \bottomrule
  \end{tabular}
\end{table}

Table~\ref{tab:vector-product-coverage} shows empirical simultaneous
coverage above the nominal level in all settings
{ for VEC-GM and VEC-WWR-LIL}. VEC-WWR-LIL covers
all 2000 paths in each setting; the Gaussian and standardized $t_5$
results are similar.

{
With the same bootstrap configuration and component kernels divided by
$\sqrt{2}A_{jj}$, Boot-AVT-Multi falls short of the nominal coverage
under both observation models in Table~\ref{tab:vector-product-coverage}.
}

\section{Concluding remarks}
\label{sec:discussion}
We develop a martingale approach to sequential inference for degenerate
degree-two kernels. Row normalization preserves the target parameter
and enables a strong Gaussian approximation on the harmonic variance
clock, at the cost of logarithmic variance inflation. Combined with
consistent variance and covariance estimation, this approximation yields
scalar confidence sequences and multi-dimensional Mahalanobis confidence
ellipsoids with asymptotic time-uniform coverage, together with one-sided
sequential tests. In the MMD and HSIC
experiments, the martingale procedures are less conservative than SAGE
and attain higher empirical power against moderate alternatives at the
same nominal level.

The results also point to two directions for further research. The analysis covers fixed kernels; empirical centering, as in practical implementations of HSIC
\citep{GrettonEtAl2005}, remains outside this scope and requires time-uniform control of the resulting
centering error. The vector theory is also restricted to multi-dimensional kernels with fixed $d$.
Extending it to high-dimensional settings may build on dimension-explicit
martingale couplings \citep{CattaneoMasiniUnderwood2025} and regularized
covariance estimation \citep{BickelLevina2008}, but their integration into
time-uniform inference remains open.

\section*{Acknowledgment}
\par  Qirui Hu's research was supported by National Natural Science Foundation of China (NSFC) (Grant Nos.12601520), the Shanghai Engineering Research Center of Finance Intelligence (Grant No.~19DZ2254600) and by TRR 391 \textit{Spatio-temporal Statistics for the Transition of Energy and Transport} (Project number 520388526) funded by the Deutsche Forschungsgemeinschaft (DFG, German Research Foundation).

\clearpage
\appendix

\section{Proofs of scalar results}
\label{sec:proofs-from-main}
\subsection{Variance clock}
\begin{proof}[Proof of Theorem~\ref{thm:scalar-clock}]
Let
\[
  q(x,y)=\E\left\{h_0(Z,x)h_0(Z,y)\right\},
  \qquad Z\sim P\text{ independent}.
\]
Then \(q\) is symmetric. By Jensen's inequality,
\[
  \E q^2(X_1,X_2)
  \le
  \E\left\{h_0^2(Z,X_1)h_0^2(Z,X_2)\right\}
  \le
  \E h_0^4(X_1,X_2)
  <\infty,
\]
and similarly \(\E q^2(X_1,X_1)<\infty\). Moreover, \(q\) is canonical because
\[
  \E\left\{q(X_1,y)\right\}
  =
  \E\left[
    h_0(Z,y)\E\left\{h_0(Z,X_1)\mid Z\right\}
  \right]
  =0.
\]
Also
\[
  \E q(X_1,X_1)
  =
  \E h_0^2(X_1,X_2)
  =
  \tau^2 .
\]

Recall that
\[
  Y_i=\frac1{i-1}\sum_{j<i}h_0(X_i,X_j),
  \qquad
  v_i=\E\left(Y_i^2\mid\calF_{i-1}\right).
\]

Since \(X_i\) is independent of \(\calF_{i-1}\) and the kernel is canonical,
\[
  \E\left(Y_i\mid\calF_{i-1}\right)
  =
  \frac1{i-1}\sum_{j<i}
  \E\left\{h_0(X_i,X_j)\mid\calF_{i-1}\right\}
  =
  \frac1{i-1}\sum_{j<i}
  \E\left\{h_0(X_i,X_j)\mid X_j\right\}
  =0.
\]
Thus \((Y_i,\calF_i)\) is a martingale-difference sequence and also
\[
\begin{aligned}
  \bracket{M}_n
  =
  \sum_{i=2}^n \E\left(Y_i^2\mid\calF_{i-1}\right)=
  \sum_{i=2}^n \frac1{(i-1)^2}
  \sum_{1\le j,k<i}q(X_j,X_k)
  =
  \sum_{j,k=1}^{n-1}a_{jk}^{(n)}q(X_j,X_k),
\end{aligned}
\]
where \(a_{jk}^{(n)}=\sum_{r=\max(j,k)}^{n-1}r^{-2}\). By the canonicity of
\(q\),
\[
  \E\bracket{M}_n
  =
  \tau^2\sum_{j=1}^{n-1}a_{jj}^{(n)}
  =
  \tau^2\sum_{r=1}^{n-1}\frac1r
  =
  \tau^2T_n .
\]

Write
\[
  \bracket{M}_n
  =
  \sum_{j=1}^{n-1}a_{jj}^{(n)}q(X_j,X_j)
  +
  2\sum_{1\le j<k\le n-1}a_{jk}^{(n)}q(X_j,X_k)
  =:A_n+B_n .
\]
For the diagonal part, \(a_{jj}^{(n)}\lesssim j^{-1}\) and
\(\E q^2(X_1,X_1)<\infty\), so $\Var(A_n)
  \lesssim
  \sum_{j=1}^{\infty}j^{-2}
  =
  O(1).$
For the off-diagonal part, the canonicity of \(q\) gives the needed
orthogonality, and
\[
\begin{aligned}
  \Var(B_n)
  &\lesssim
  \sum_{1\le j<k\le n-1}
  \left\{a_{jk}^{(n)}\right\}^2
  \lesssim
  \sum_{k=2}^{n-1}\sum_{j=1}^{k-1}\frac1{k^2}
  =
  O(T_n).
\end{aligned}
\]
For \(j<k\),
\[
  \E\left\{q(X_j,X_j)q(X_j,X_k)\right\}
  =
  \E\left[
    q(X_j,X_j)\E\left\{q(X_j,X_k)\mid X_j\right\}
  \right]
  =0
\]
by canonicity of \(q\). The same conditioning argument eliminates the other
diagonal--off-diagonal cross terms, so
\[
  \Cov(A_n,B_n)=0.
\]
Consequently, with \(\Delta_n=\bracket{M}_n-\tau^2T_n\),
\[
  \E\Delta_n=0,
  \qquad
  \Var(\Delta_n)=\Var(A_n)+\Var(B_n)=O(T_n).
\]

Choose a sparse grid \(n_k=\lfloor\exp(k^p)\rfloor\) with \(p>2\). Since
\(T_n\sim\log n\), \(T_{n_k}\asymp k^p\) and
\(\log T_{n_k}\asymp\log k\). Chebyshev's inequality and Borel--Cantelli give
\[
  \Delta_{n_k}
  =
  o_{a.s.}\left(\frac{T_{n_k}}{\log T_{n_k}}\right).
\]
If \(n_k\le n\le n_{k+1}\), then monotonicity gives
\[
  \Delta_{n_k}-\tau^2\left(T_n-T_{n_k}\right)
  \le
  \Delta_n
  \le
  \Delta_{n_{k+1}}+\tau^2\left(T_{n_{k+1}}-T_n\right).
\]
After division by \(T_n/\log T_n\), the grid errors vanish by the preceding
display and the deterministic increments vanish because
\[
  T_{n_{k+1}}-T_{n_k}=O(k^{p-1})
  =
  o\left(\frac{k^p}{\log k}\right)
  =
  o\left(\frac{T_{n_k}}{\log T_{n_k}}\right).
\]
It follows that
\[
  \bracket{M}_n-\tau^2T_n
  =
  o_{a.s.}\left(\frac{T_n}{\log T_n}\right),
  \qquad
  \frac{\bracket{M}_n}{T_n}\asto\tau^2 .
\]

For the optional quadratic variation, define
\[
  D_i=Y_i^2-\E\left(Y_i^2\mid\calF_{i-1}\right),\qquad i\ge2.
\]
Then \((D_i,\calF_i)\) is a martingale-difference sequence. Conditionally on
\(X_i\), the variables $h_0(X_i,X_j)$, $1\le j<i$, are
independent and centered. The fourth-moment formula for sums of independent
centered random variables and Jensen's inequality give $\E Y_i^4=O(i^{-2}).$ Therefore $\sum_{i=2}^{\infty}\E D_i^2<\infty,$ so \(\sum_{i=2}^nD_i\) converges almost surely and in \(L^2\). Hence
\[
  \qv{M}_n-\bracket{M}_n
  =
  O_{a.s.}(1)
  =
  o_{a.s.}\left(\frac{T_n}{\log T_n}\right),
\]
and
\[
  \frac{\qv{M}_n}{T_n}\asto\tau^2 .
\]

Finally, since \(M_n=\sum_{i=2}^nY_i=(n-1)(\widetilde U_n-\theta)\),
\[
\begin{aligned}
\sum_{i=2}^n
\left\{
\frac1{i-1}\sum_{j<i}h(X_i,X_j)-\widetilde U_n
\right\}^2
&=
\sum_{i=2}^n
\left(
Y_i-\frac{M_n}{n-1}
\right)^2\\
&=
\qv{M}_n-\frac{M_n^2}{n-1}.
\end{aligned}
\]
Thus
\[
\widehat\tau_n^2
=
\frac{\qv{M}_n}{T_n}
-
\frac{M_n^2}{(n-1)T_n}.
\]
Since
\[
  \sum_{i=2}^{\infty}\E\left(\frac{Y_i}{\sqrt{i}}\right)^2
  =
  \tau^2\sum_{i=2}^{\infty}\frac1{i(i-1)}
  <\infty,
\]
the martingale \(L_n=\sum_{i=2}^nY_i/\sqrt{i}\) is \(L^2\)-bounded. The
martingale convergence theorem gives \(L_n\to L_\infty\) almost surely.
Kronecker's lemma, applied with the increasing weights \(\sqrt{i}\), yields
\[
  \frac1{\sqrt n}\sum_{i=2}^nY_i
  =
  \frac{M_n}{\sqrt n}
  \to0
  \qquad\text{a.s.}
\]
Hence
\[
  \frac{M_n^2}{(n-1)T_n}=o_{a.s.}(1),
\]
and \(\widehat\tau_n^2\asto\tau^2\).
\end{proof}

\subsection{Strong approximation}

The following embedding is given in \citet[Appendix~I]{HallHeyde1980}.

\begin{lemma}[Martingale Skorokhod embedding]
\label{lem:HH-skorokhod}
Let $S_0=0$ and $S_n=\sum_{i=1}^n X_i$ be a square-integrable
martingale with respect to $\{\mathcal F_n\}$. There exist a standard
Brownian motion $W$, an increasing sequence of embedding times
$0=R_0\le R_1\le\cdots$, and a filtration $(\mathcal G_n)$ such that
$S'_n=W(R_n)$ is a martingale with respect to $(\mathcal G_n)$ and
$(S'_n)_{n\ge0}$ has the same law as $(S_n)_{n\ge0}$.
Writing $\rho_i=R_i-R_{i-1}$ and $X'_i=S'_i-S'_{i-1}$, we have
\[
\mathbb E(\rho_i\mid\mathcal G_{i-1})
 =\mathbb E\bigl((X'_i)^2\mid\mathcal G_{i-1}\bigr),
\qquad
\mathbb E(\rho_i^2\mid\mathcal G_{i-1})
 \le C\,\mathbb E\bigl(|X'_i|^4\mid\mathcal G_{i-1}\bigr)
\]
for a universal constant $C$.
\end{lemma}

\begin{lemma}[Random-time consistency]\label{lem:random-time-consistency}
\leavevmode\newline
Under the canonical condition and \(E|h_0(X_1,X_2)|^4<\infty\), the embedding
times for \(M_n=\sum_{i=2}^nY_i\) satisfy
\[
  \left|R_n-\tau^2T_n\right|
  =
  o_{a.s.}\left(\frac{T_n}{\log T_n}\right).
\]
In particular, \(R_n/T_n\asto\tau^2\).
\end{lemma}

\begin{proof}
Apply Lemma~\ref{lem:HH-skorokhod} to the martingale differences \(Y_i\). Let
\(M_n'=B(R_n)\) be the embedded martingale copy, where
\(R_n=\sum_{i=2}^n\rho_i\), and set \(Y_i'=M_i'-M_{i-1}'\). In the embedded
filtration \((\calG_i)\), define
\[
  v_i'=\E\left\{(Y_i')^2\mid\calG_{i-1}\right\},
  \qquad
  \bracket{M'}_n=\sum_{i=2}^n v_i' .
\]
The equality in law of \((Y_i')\) and \((Y_i)\) gives
\[
  \sum_{i=2}^{\infty}\E\left\{(Y_i')^4\right\}
  =
  \sum_{i=2}^{\infty}\E Y_i^4
  <\infty.
\]
Write \(D_i'=(Y_i')^2-v_i'\). Then \((D_i',\calG_i)\) is a
martingale-difference sequence and
\[
  \sum_{i=2}^{\infty}\E(D_i')^2
  \le
  \sum_{i=2}^{\infty}\E\left\{(Y_i')^4\right\}
  <\infty.
\]
Thus
\[
  \qv{M'}_n-\bracket{M'}_n
  =
  \sum_{i=2}^nD_i'
  =O_{a.s.}(1).
\]
The proof of Theorem~\ref{thm:scalar-clock} establishes
\[
  \qv{M}_n-\tau^2T_n
  =
  o_{a.s.}\left(\frac{T_n}{\log T_n}\right).
\]
Because \((Y_i')\) and \((Y_i)\) have the same finite-dimensional
distributions, this optional-quadratic-variation rate transfers to the
embedded copy:
\[
  \qv{M'}_n-\tau^2T_n
  =
  o_{a.s.}\left(\frac{T_n}{\log T_n}\right).
\]
Consequently,
\[
  \bracket{M'}_n-\tau^2T_n
  =
  o_{a.s.}\left(\frac{T_n}{\log T_n}\right).
\]

Moreover,
\[
  \E\left(\rho_i\mid\calG_{i-1}\right)=v_i',
\]
so \(\sum_{i=2}^n(\rho_i-v_i')\) is a martingale. The conditional
second-moment bound in the embedding lemma gives
\[
  \sum_{i=2}^{\infty}\E(\rho_i-v_i')^2
  \le
  C\sum_{i=2}^{\infty}\E (Y_i')^4
  <\infty.
\]
Hence
\[
  R_n-\bracket{M'}_n=O_{a.s.}(1).
\]
Combining the last two displays gives the asserted rate.
\end{proof}

\begin{proof}[Proof of Theorem~\ref{thm:scalar-asip}]
By the Hall--Heyde Skorokhod embedding applied to the martingale differences
\(Y_i\), we may work on an enlarged probability space on which $M_n=B(R_n)$ where \(B\) is a standard Brownian motion and \(R_n\) is the embedding clock.
By Lemma~\ref{lem:random-time-consistency},
\[
    |R_n-\tau^2T_n|
    =
    o_{a.s.}\left(\frac{T_n}{\log T_n}\right).
\]
Let
\[
    \omega_B(T,h)
    =
    \sup_{\substack{0\le s,t\le T\\ |s-t|\le h}}
    |B(s)-B(t)|.
\]
The Brownian modulus of continuity gives
\[
    \omega_B(T,h)
    =
    O_{a.s.}\left(
    \sqrt{h\left\{\log(T/h)+\log\log T\right\}}
    \right)
\]
uniformly as \(T\to\infty\), for \(0<h\le T/2\). Since
\(R_n/(\tau^2T_n)\to1\) almost surely, with
\(h_n=|R_n-\tau^2T_n|\),
choose a deterministic constant \(C>\tau^2\). Eventually,
\[
  R_n\vee(\tau^2T_n)\le CT_n,
  \qquad
  h_n=o_{a.s.}\left(\frac{T_n}{\log T_n}\right).
\]
Thus \(h_n\le CT_n/2\) eventually, and the modulus bound, applied with
\(T=CT_n\), gives
\[
  h_n\left\{\log(CT_n/h_n)+\log\log(CT_n)\right\}
  =
  o_{a.s.}(T_n).
\]
Here the display is interpreted trivially when \(h_n=0\). Therefore
\[
    |B(R_n)-B(\tau^2T_n)|
    =
    o_{a.s.}(\sqrt{T_n}).
\]
Defining \(W(t)=\tau^{-1}B(\tau^2t)\), we obtain a standard Brownian motion
\(W\). This is the standard martingale strong-approximation argument obtained
from Skorokhod embedding plus Brownian modulus of continuity; see
\citet[Appendix I]{HallHeyde1980}. Hence
\[
    M_n
    =
    \tau W(T_n)+o_{a.s.}(\sqrt{T_n}).
\]
Division by \(n-1\) gives the equivalent display for
\(\widetilde U_n-\theta\).
\end{proof}

\subsection{Time-uniform coverage}

\begin{lemma}[Scalar Brownian calibration]
\label{lem:scalar-brownian-calibration}
Both \eqref{eq:gm-boundary} and \eqref{eq:lil-boundary} are two-sided
\((1-\alpha)\) Brownian motion boundaries in the sense of
Definition~\ref{def:scalar-boundary}, for every \(\alpha\in(0,1)\) and
\(t_0>0\). The LIL formula permits every fixed \(\eta>1\) and \(s>1\).
\end{lemma}

\begin{proof}
For GM, let \(c=g^{-1}(\alpha)\). On \(t\ge t_0\), the process $L_t
  =
  t^{-1/2}
  \exp\left\{{W(t)^2}/{2t}\right\}$ is a nonnegative martingale conditional on \(W(t_0)\).

A closed crossing of the GM boundary is equivalent to $L_t
  \ge
  t_0^{-1/2}e^{c^2/2}$ for some \(t\ge t_0\). With \(Z=W(t_0)/\sqrt{t_0}\), the conditional
maximal inequality and then Gaussian integration give
\begin{align*}
    \Pp\left\{\exists t\ge t_0:|W(t)|\ge u_{\alpha,t_0}^{\rm GM}(t)\right\}=&\E\left[\Pp\left\{
    \sup_{t\ge t_0}L_t\ge
    t_0^{-1/2}e^{c^2/2}\mid W(t_0)
  \right\}\right]\\&\le
  \E\min\left\{
    1,e^{(Z^2-c^2)/2}
  \right\}
  =
  2\{1-\Phi(c)+c\phi(c)\}
  =
  \alpha.
\end{align*}

This is the classical mixture-martingale boundary-crossing argument of \cite{Robbins1970b}.

For LIL, Brownian scaling reduces the proof to \(t_0=1\).
For \(k\ge0\), set $A_k
  =
  \log\left\{
{2\zeta(s)(k+1)^s}/{\alpha}
  \right\}$ and $\lambda_k
  =
  \sqrt{
    {2A_k}/{\eta^{k+1/2}}
  }.$
Consider the exponential martingale $\exp\left\{
    \lambda_k W(t)
    -{\lambda_k^2t}/{2}
  \right\}$
and let $L_k(t)
  :={A_k}/{\lambda_k}
  +{\lambda_k t}/{2},$ then
by Ville's inequality,
\[
  \Pp\left\{
    \exists t\ge0:
    W(t)\ge L_k(t)
  \right\}
  \le
  \Pp\left\{
    \sup_{t\ge0}
    \exp\left(
      \lambda_kW(t)
      -
      \frac{\lambda_k^2t}{2}
    \right)
    \ge e^{A_k}
  \right\}
  \le
  e^{-A_k}.
\]

For \(t\in[\eta^k,\eta^{k+1}]\),
\[
  \frac{L_k(t)}{\sqrt{A_kt}}
  =
  \frac{
    \sqrt{\eta^{k+1/2}/t}
    +
    \sqrt{t/\eta^{k+1/2}}
  }{\sqrt2}
  \le
  \frac{
    \eta^{1/4}
    +
    \eta^{-1/4}
  }{\sqrt2}.
\]
Moreover, the expression in square brackets in
\eqref{eq:lil-boundary} is at least \(A_k\) on this epoch. Indeed,
since \(t\ge\eta^k\),
\[
\begin{aligned}
  s\log\log\{\max(\eta t,e)\}
  +
  \log\frac{2\zeta(s)}
  {\alpha(\log\eta)^s}&\ge
  s\log\{(k+1)\log\eta\}
  +
  \log\frac{2\zeta(s)}
  {\alpha(\log\eta)^s}
\\
  &=
  \log\left\{
    \frac{2\zeta(s)(k+1)^s}{\alpha}
  \right\}
  =
  A_k.
\end{aligned}
\]
Thus, on the \(k\)-th epoch, $u_{\alpha,1}^{\rm LIL}(t)
  \ge
  L_k(t).$
Therefore a crossing of the upper LIL boundary during this epoch implies
a crossing of \(L_k\), and
\[
  \Pp\left\{
    \exists t\in[\eta^k,\eta^{k+1}]:
    W(t)\ge
    u_{\alpha,1}^{\rm LIL}(t)
  \right\}
  \le
  e^{-A_k}.
\]
Summing over \(k\ge0\),
\[
\begin{aligned}
  \Pp\left\{
    \exists t\ge1:
    W(t)\ge
    u_{\alpha,1}^{\rm LIL}(t)
  \right\}
  \le
  \sum_{k=0}^\infty e^{-A_k} =
  \frac{\alpha}{2\zeta(s)}
  \sum_{k=0}^\infty\frac1{(k+1)^s}
  =
  \frac{\alpha}{2}.
\end{aligned}
\]
By symmetry, the same bound holds for crossing the lower boundary.
Hence
\[
  \Pp\left\{
    \exists t\ge1:
    |W(t)|
    \ge
    u_{\alpha,1}^{\rm LIL}(t)
  \right\}
  \le
  \alpha.
\]
Brownian scaling gives the same conclusion for every \(t_0>0\).
\end{proof}

\begin{proof}[Proof of Theorem~\ref{thm:scalar-asympcs}]
By Theorem~\ref{thm:scalar-asip}, we have $M_n
  =
  \tau W(T_n)+r_n=\tau W(T_n)+o_{a.s.}(\sqrt{T_n})$.

For either boundary, $u_{\alpha,t_0}(t_0r)
  =
  \sqrt{t_0}\,
  u_{\alpha,1}(r)$ for $r\ge1$ and $c_0
  :=
  \inf_{r\ge1}
  {u_{\alpha,1}(r)}/{\sqrt r}
  >0.$
Hence
\[
  u_{\alpha,T_m}(T_n)
  \ge
  c_0\sqrt{T_n},
  \qquad n\ge m.
\]

Define \[
  \Delta_m
  =
  \sup_{n\ge m}
  \left[
    \frac{|r_n|}
    {\tau u_{\alpha,T_m}(T_n)}
    +
    \left|
      \frac{\widehat\tau_n}{\tau}-1
    \right|
  \right].
\]
Since \(\widehat\tau_n/\tau\asto1\), we have
\[
  \Delta_m
  \le
  \frac1{\tau c_0}
  \sup_{n\ge m}
  \frac{|r_n|}{\sqrt{T_n}}
  +
  \sup_{n\ge m}
  \left|
    \frac{\widehat\tau_n}{\tau}-1
  \right|
  \asto0.
\]

Consider the feasible and ideal crossing events $  E_m
  =
  \left\{
    \exists n\ge m:
    |M_n|
    >
    \widehat\tau_n
    u_{\alpha,T_m}(T_n)
  \right\}$
and $A_m
  =
  \left\{
    \exists t\ge T_m:
    |W(t)|
    \ge
    u_{\alpha,T_m}(t)
  \right\}.$
By Lemma~\ref{lem:scalar-brownian-calibration}, $\Pp(A_m)\le\alpha$ for every \(m\).

If \(E_m\) occurs, then for some \(n\ge m\),
\[
\begin{aligned}
  |W(T_n)|
  &\ge
  \frac{|M_n|-|r_n|}{\tau} >
  \frac{\widehat\tau_n}{\tau}
  u_{\alpha,T_m}(T_n)
  -
  \frac{|r_n|}{\tau} \\
  &\ge
  (1-\Delta_m)
  u_{\alpha,T_m}(T_n).
\end{aligned}
\]
Thus, for any fixed \(\delta\in(0,1)\),
\[
\begin{aligned}
  \Pp(E_m\setminus A_m)
  \le
  \Pp(\Delta_m>\delta) +
  \Pp\left\{
    A_m^c,\ %
    \exists t\ge T_m:
    |W(t)|
    >
    (1-\delta)
    u_{\alpha,T_m}(t)
  \right\}.
\end{aligned}
\]
By Brownian scaling and the scaling identity of the boundary, the second
probability equals
\[
  \Pp\left\{
    \forall r\ge1:
    |W(r)|<u_{\alpha,1}(r),\
    \exists r\ge1:
    |W(r)|
    >
    (1-\delta)u_{\alpha,1}(r)
  \right\}.
\]

We next show that this probability vanishes as \(\delta\downarrow0\).
For the Gaussian-mixture boundary, ${
    u_{\alpha,1}^{\rm GM}(r)
  }/{\sqrt{2r\log\log r}}
  \rightarrow\infty,$
whereas for the stitched-LIL boundary, ${
    \{u_{\alpha,1}^{\rm LIL}(r)\}^2
  }/{
    2r\log\log r
  }
  \rightarrow{
    s(\eta^{1/4}+\eta^{-1/4})^2
  }/{4}
  >1.$
Therefore, by the Brownian law of the iterated logarithm, we have
\[
  \limsup_{r\to\infty}
  \frac{|W(r)|}{u_{\alpha,1}(r)}
  <1
  \qquad\text{a.s.}
\]
On the event $\left\{
    |W(r)|<u_{\alpha,1}(r)
    \text{ for all } r\ge1
  \right\},$ the preceding limsup bound implies that there exist
\(\varepsilon>0\) and \(R<\infty\) such that ${|W(r)|}/{u_{\alpha,1}(r)}
  \le
  1-\varepsilon$ for $r\ge R$.
On the compact interval \([1,R]\), continuity of \(W\) and
\(u_{\alpha,1}\) gives $\sup_{1\le r\le R}
  {|W(r)|}/{u_{\alpha,1}(r)}
  <1,$
since the ratio is strictly below \(1\) at every point. Therefore on this event
\[
  \sup_{r\ge1}
  \frac{|W(r)|}{u_{\alpha,1}(r)}
  <1
  \qquad\text{a.s..}
\]
Consequently, as $\delta\downarrow0$ we have
\[
\Pp\left\{
  \forall r\ge1:
  |W(r)|<u_{\alpha,1}(r),\
  \exists r\ge1:
  |W(r)|
  >
  (1-\delta)u_{\alpha,1}(r)
\right\}
\rightarrow0.
\]
Since \(\Delta_m\asto0\), it follows that for every fixed
\(\delta\in(0,1)\),
\[
\begin{aligned}
  \limsup_{m\to\infty}\Pp(E_m)
  &\le
  \alpha
  +
  \Pp\left\{
    \forall r\ge1:
    |W(r)|<u_{\alpha,1}(r),\
    \exists r\ge1:
    |W(r)|
    >
    (1-\delta)u_{\alpha,1}(r)
  \right\}.
\end{aligned}
\]
Letting \(\delta\downarrow0\) gives $\limsup_{m\to\infty}\Pp(E_m)
  \le\alpha.$

This is exactly
\[
  \liminf_{m\to\infty}
  \Pp\left\{
    \forall n\ge m:
    \theta
    \in
    C_{m,n}^{\rm MG}(\alpha)
  \right\}
  \ge
  1-\alpha,
\]
which proves the theorem.
\end{proof}

\subsection{Asymptotic confidence sequences}

\begin{proof}[Proof of Remark~\ref{rem:scalar-asympcs}]
Fix $m\ge2$. On the coupling in Theorem~\ref{thm:scalar-asip}, as $M_n=\tau W(T_n)+r_n$ where $r_n=o_{a.s.}(\sqrt{T_n})$, we have
\[
  A_m:=\sup_{n\ge m}\frac{|r_n|}{\sqrt{T_n}}<\infty
  \qquad\text{a.s.}.
\]
Fix any $\delta\in(0,\alpha)$. Since $A_m<\infty$ almost surely, there
exists a deterministic $\varepsilon>0$ such that
\[
  \Pp(A_m>\varepsilon)\le\delta.
\]
Put $a=\alpha-\delta$ and define
\[
  w^*_{m,n}
  =
  \frac{\tau u_{a,T_m}(T_n)+\varepsilon\sqrt{T_n}}{n-1},
  \qquad n\ge m,
\]
and $C^*_{m,n}(\alpha)
  =
  [\,\widetilde U_n\pm w^*_{m,n}\,].$ On the event $ \{A_m\le\varepsilon\}
  \cap
  \bigl\{
    |W(T_n)|\le u_{a,T_m}(T_n)
    \text{ for all }n\ge m
  \bigr\},$ we have for every $n\ge m$,
\[
  |M_n|
  \le
  \tau u_{a,T_m}(T_n)+\varepsilon\sqrt{T_n},
\]
and therefore $\theta\in C^*_{m,n}(\alpha)$ for all $n\ge m$.
By Lemma~\ref{lem:scalar-brownian-calibration} and a union bound,
\[
  \Pp\{\theta\in C^*_{m,n}(\alpha)
  \text{ for all }n\ge m\}
  \ge
  1-a-\delta
  =
  1-\alpha.
\]
Thus $\{C^*_{m,n}(\alpha):n\ge m\}$ is a nonasymptotic
$(1-\alpha)$-confidence sequence. The constants $\tau$ and
$\varepsilon$ may depend on the underlying law and need not be available
for implementation.

It remains to compare its radius with that of the feasible sequence. Recall that $w_{m,n}
  ={\widehat\tau_n u_{\alpha,T_m}(T_n)}/{(n-1)},$
we have
\[
  \frac{w^*_{m,n}}{w_{m,n}}
  =
  \frac{\tau}{\widehat\tau_n}
  \frac{u_{a,T_m}(T_n)}{u_{\alpha,T_m}(T_n)}
  +
  \frac{\varepsilon}{\widehat\tau_n}
  \frac{\sqrt{T_n}}{u_{\alpha,T_m}(T_n)}
\]
almost surely for all sufficiently large $n$.
By Theorem~\ref{thm:scalar-clock},
$\widehat\tau_n\asto\tau>0$. For either the Gaussian-mixture or
stitched-LIL boundary, with $m$ fixed, ${u_{a,T_m}(T_n)}/{u_{\alpha,T_m}(T_n)}
  \rightarrow1$ and ${\sqrt{T_n}}/{u_{\alpha,T_m}(T_n)}
  \rightarrow0$. Consequently,
\[
  \frac{w^*_{m,n}}{w_{m,n}}
  \asto 1,
  \qquad n\to\infty.
\]
Hence $\{C^{\mathrm{MG}}_{m,n}(\alpha):n\ge m\}$ satisfies Definition~\ref{def:comparison}.
\end{proof}

\subsection{One-sided sequential testing}
\begin{lemma}[One-sided Brownian boundaries]
\label{lem:scalar-one-sided-boundaries}
Let $t_0>0$. The boundaries in Remark~\ref{rem:one-sided-testing} satisfy
\[
  \Pp\left\{\exists t\ge t_0:
    W(t)\ge u_{\alpha,t_0}^+(t)\right\}\le\alpha,
\]
where $W$ is a standard Brownian motion. For the LIL boundary, this holds
for every $\alpha\in(0,1)$, $\eta>1$, and $s>1$; for the GM boundary, it
holds for $0<\alpha<g_+(0)=1/2+1/\pi$.
\end{lemma}

\begin{proof}
For the GM boundary, write $c=g_+^{-1}(\alpha)>0$. The function $g_+$ is
continuous and strictly decreasing on $[0,\infty)$, since
\[
  g_+'(a)
  =-\phi(a)\left\{a+\frac{\phi(a)}{\Phi(a)}\right\}^2<0,
\]
and $g_+(a)\to0$ as $a\to\infty$. Hence $g_+(c)=\alpha$.
For $t\ge t_0$, define
\[
  L_t^+
  =t^{-1/2}\exp\left\{\frac{W(t)^2}{2t}\right\}
    \Phi\left(\frac{W(t)}{\sqrt t}\right)
  =\frac1{\sqrt{2\pi}}\int_0^\infty
    \exp\left\{\lambda W(t)-\frac{\lambda^2t}{2}\right\}\,d\lambda.
\]
Conditionally on $W(t_0)$, this is a nonnegative martingale on
$[t_0,\infty)$. A crossing of the upper GM boundary implies
\[
  L_t^+\ge t_0^{-1/2}e^{c^2/2}\Phi(c)
\]
at the crossing time. With $Z=W(t_0)/\sqrt{t_0}\sim N(0,1)$, the
conditional maximal inequality therefore yields
\begin{align*}
  \Pp\left\{\exists t\ge t_0:
    W(t)\ge u_{\alpha,t_0}^{+,\mathrm{GM}}(t)\right\}
  &\le
  \E\min\left\{1,
    e^{(Z^2-c^2)/2}\frac{\Phi(Z)}{\Phi(c)}\right\} \\
  &=1-\Phi(c)
    +\frac{\phi(c)}{\Phi(c)}\int_{-\infty}^{c}\Phi(z)\,dz \\
  &=1-\Phi(c)+c\phi(c)+\frac{\phi(c)^2}{\Phi(c)}
   =\alpha.
\end{align*}
The second line uses the strict monotonicity of
$z\mapsto e^{z^2/2}\Phi(z)$, which follows from its positive-integral
representation, and the third uses
$\int_{-\infty}^{c}\Phi(z)\,dz=c\Phi(c)+\phi(c)$.

For the LIL boundary, Brownian scaling reduces the proof to $t_0=1$.
For $k\ge0$, put
\[
  A_k=\log\{\zeta(s)(k+1)^s/\alpha\},\qquad
  \lambda_k=\sqrt{2A_k/\eta^{k+1/2}},\qquad
  L_k(t)=A_k/\lambda_k+\lambda_k t/2.
\]
Ville's inequality for
$\exp\{\lambda_kW(t)-\lambda_k^2t/2\}$ gives
\[
  \Pp\{\exists t\ge0:W(t)\ge L_k(t)\}\le e^{-A_k}.
\]
For $t\in[\eta^k,\eta^{k+1}]$,
\[
  \frac{L_k(t)}{\sqrt{A_kt}}
  =\frac{\sqrt{\eta^{k+1/2}/t}+\sqrt{t/\eta^{k+1/2}}}{\sqrt2}
  \le\frac{\eta^{1/4}+\eta^{-1/4}}{\sqrt2}.
\]
Moreover,
\[
  s\log\log\{\max(\eta t,e)\}
  +\log\frac{\zeta(s)}{\alpha(\log\eta)^s}
  \ge\log\frac{\zeta(s)(k+1)^s}{\alpha}=A_k.
\]
Thus $u_{\alpha,1}^{+,\mathrm{LIL}}(t)\ge L_k(t)$ on this epoch.
Summing the epochwise crossing probabilities gives
\[
  \Pp\left\{\exists t\ge1:
    W(t)\ge u_{\alpha,1}^{+,\mathrm{LIL}}(t)\right\}
  \le\sum_{k=0}^\infty e^{-A_k}
  =\frac{\alpha}{\zeta(s)}\sum_{k=0}^\infty(k+1)^{-s}
  =\alpha.
\]
Brownian scaling gives the conclusion for arbitrary $t_0>0$.
\end{proof}

\begin{proof}[Proof of Corollary~\ref{cor:one-sided-testing}]
Work under the fixed null distribution in the statement. By
Theorem~\ref{thm:scalar-asip},
\[
  M_n=(n-1)(\widetilde U_n-\theta_0)
  =\tau W(T_n)+r_n,\qquad r_n=o_{a.s.}(\sqrt{T_n}),
\]
and Theorem~\ref{thm:scalar-clock} gives $\widehat\tau_n/\tau\asto1$.
For either one-sided boundary,
\[
  u_{\alpha,t_0}^+(t_0r)=\sqrt{t_0}\,u_{\alpha,1}^+(r),\qquad
  c_0:=\inf_{r\ge1}\frac{u_{\alpha,1}^+(r)}{\sqrt r}>0.
\]
Consequently,
\[
  \Delta_m
  :=\sup_{n\ge m}\left\{
      \frac{|r_n|}{\tau u_{\alpha,T_m}^+(T_n)}
      +\left|\frac{\widehat\tau_n}{\tau}-1\right|
    \right\}
  \asto0.
\]
If the test rejects at some $n\ge m$, then
$W(T_n)>(1-\Delta_m)u_{\alpha,T_m}^+(T_n)$.
Lemma~\ref{lem:scalar-one-sided-boundaries} bounds the ideal closed-crossing
probability by $\alpha$. To transfer this bound, note that
\[
  \frac{u_{\alpha,1}^{+,\mathrm{GM}}(r)}{\sqrt{2r\log\log r}}
  \longrightarrow\infty,\qquad
  \frac{\{u_{\alpha,1}^{+,\mathrm{LIL}}(r)\}^2}{2r\log\log r}
  \longrightarrow\frac{s(\eta^{1/4}+\eta^{-1/4})^2}{4}>1.
\]
The Brownian LIL therefore gives
$\limsup_{r\to\infty}W(r)/u_{\alpha,1}^+(r)<1$ almost surely.
On the event of no ideal closed crossing, continuity on compact intervals
then implies $\sup_{r\ge1}W(r)/u_{\alpha,1}^+(r)<1$ almost surely.
The same scaling and vanishing-perturbation argument as in the proof of
Theorem~\ref{thm:scalar-asympcs}, with $W$ in place of $|W|$, now yields
the asserted upper bound on the limiting rejection probability.
\end{proof}

\subsection{Fixed-alternative consistency}
\label{app:one-sided-fixed-alternative}

\begin{proof}[Proof of Proposition~\ref{prop:one-sided-fixed-alternative}]
Define $\mu(x)=\int h(x,y)\,P(dy),$ $a_i=\sum_{j<i}h(X_i,X_j)/(i-1)$ and $e_i=a_i-\mu(X_i)$ for $i\ge2.$ Conditional Jensen's inequality gives
$\E|\mu(X_1)-\theta|^4\le\E|h_0(X_1,X_2)|^4<\infty$.
Writing $\xi=h(X_1,X_2)-\mu(X_1)$, we also have
\[
  K_4:=\E|\xi|^4
  \le16\E|h_0(X_1,X_2)|^4<\infty.
\]
Conditionally on $X_i$, the variables
$h(X_i,X_j)-\mu(X_i)$, $j<i$, are independent and centered.
Expanding the conditional fourth moment and then taking expectations
yields
\begin{align*}
  \E|e_i|^4
  =\frac{K_4}{(i-1)^3}
    +\frac{3(i-2)}{(i-1)^3}
      \E\!\left[\{\E(\xi^2\mid X_1)\}^{2}\right] \le\frac{3K_4}{(i-1)^2},
\end{align*}
where the inequality follows from conditional Jensen's inequality.
Thus, for every $\varepsilon>0$,
\[
  \sum_{i=2}^{\infty}\Pp(|e_i|>\varepsilon)
  \le\frac{3K_4}{\varepsilon^4}
       \sum_{i=2}^{\infty}\frac1{(i-1)^2}
  <\infty.
\]
The Borel--Cantelli lemma implies $e_i\asto0$;
independence between rows is not needed.
The strong law for $\mu(X_i)$ and Ces\`aro averaging therefore give
\[
  \widetilde U_n
  =\frac1{n-1}\sum_{i=2}^n\mu(X_i)
   +\frac1{n-1}\sum_{i=2}^n e_i
  \asto\theta.
\]

To control the scale estimator, note that $\sum_{i=2}^n e_i^2/(n-1)\asto0$ and $\sum_{i=2}^n\mu(X_i)^2/(n-1)
  \asto\E\{\mu(X_1)^2\}<\infty$ by Ces\`aro averaging and the strong law, respectively.
Since $a_i=\mu(X_i)+e_i$, the Cauchy--Schwarz inequality yields
\begin{align*}
  \frac1{n-1}\sum_{i=2}^n
    |a_i^2-\mu(X_i)^2|
  &\le
    2\left\{\frac1{n-1}\sum_{i=2}^n\mu(X_i)^2\right\}^{1/2}
     \left\{\frac1{n-1}\sum_{i=2}^n e_i^2\right\}^{1/2} \\
  &\quad+\frac1{n-1}\sum_{i=2}^n e_i^2
  \asto0.
\end{align*}
Consequently,
\[
  \frac{T_n}{n-1}\widehat\tau_n^2
  =\frac1{n-1}\sum_{i=2}^n(a_i-\widetilde U_n)^2
  =\frac1{n-1}\sum_{i=2}^n a_i^2-\widetilde U_n^2
  \asto\Var\{\mu(X_1)\}.
\]
In particular,
$\widehat\tau_n^2=O_{\mathrm{a.s.}}(n/T_n)$,
including when $\Var\{\mu(X_1)\}=0$.

For fixed $m$ and fixed boundary parameters, the formulas in
Remark~\ref{rem:one-sided-testing} imply
\[
  \frac{\{u_{\alpha,T_m}^+(T_n)\}^2}{T_n}
  =\begin{cases}
    O(\log T_n),&\mathrm{GM},\\
    O(\log\log T_n),&\mathrm{LIL}.
  \end{cases}
\]
Since $T_n\sim\log n$, it follows that
\begin{align*}
  \frac{\widehat\tau_n u_{\alpha,T_m}^+(T_n)}{n-1}
  &=\left(\frac{T_n\widehat\tau_n^2}{n-1}\right)^{1/2}
    \frac{u_{\alpha,T_m}^+(T_n)}{\sqrt{(n-1)T_n}} \\
  &=\begin{cases}
    O_{\mathrm{a.s.}}\!\left(\sqrt{\dfrac{\log\log n}{n}}\right),
      &\mathrm{GM},\\[3pt]
    O_{\mathrm{a.s.}}\!\left(\sqrt{\dfrac{\log\log\log n}{n}}\right),
      &\mathrm{LIL},
  \end{cases}
  \asto0.
\end{align*}
Combining this with $\widetilde U_n\asto\theta$ gives
\[
  \widetilde U_n-\theta_0
  -\frac{\widehat\tau_n u_{\alpha,T_m}^+(T_n)}{n-1}
  \asto\theta-\theta_0>0.
\]
Hence the rejection inequality holds at every sufficiently large
monitoring time, proving the assertion.
The argument also covers times at which $\widehat\tau_n=0$,
since it never divides by the estimated scale.
\end{proof}

\section{Proofs of vector results}
\subsection{Covariance clock}

\begin{lemma}[Vector covariance clock and plug-in consistency]
\label{lem:vector-clock}
Under (A2),
\[
  \norm{
    \bracket{\bm M}_n-T_n\bm\Omega
  }_{\op}
  =
  O_{a.s.}(T_n^{3/4})
  =
  o_{a.s.}\left(\frac{T_n}{\log T_n}\right),
\]
where $\bracket{\bm M}_n:=\sum_{i=2}^n\E\left(
    \bm Y_i\bm Y_i\trans
    \mid\calF_{i-1}
  \right)$ and hence ${\bracket{\bm M}_n}/T_n\asto\bm\Omega$ in the operator norm. Moreover, the plug-in estimator is strongly consistent $\widehat{\bm\Omega}_n\asto\bm\Omega$.
\end{lemma}

\begin{proof}
Conditionally on
\(X_i\), the vectors \(\bm h_0(X_i,X_j)\), \(j<i\), are independent.
They are also conditionally centered, since symmetry and canonicity give
\(\E\{\bm h_0(X_i,X_j)\mid X_i\}=\bm0\). Thus, for \(P\)-a.e. fixed \(x\), applying Rosenthal inequality with exponent four coordinatewise and
using
\(\norm{\bm u}_2^4\le d\sum_{r=1}^d|u_r|^4\), we obtain
\[
\begin{aligned}
  \E\left\|
    \sum_{j<i}\bm h_0(x,X_j)
  \right\|_2^4
  \le C_d\Big[
    &(i-1)\E\norm{\bm h_0(x,X_1)}_2^4\\
    &+
    (i-1)^2
    \{\E\norm{\bm h_0(x,X_1)}_2^2\}^2
  \Big].
\end{aligned}
\]
By Jensen's inequality,
\(\{\E\norm{\bm h_0(x,X_1)}_2^2\}^2
\le\E\norm{\bm h_0(x,X_1)}_2^4\). Evaluating the preceding bound at
\(x=X_i\) and using independence of \(X_i\) and
\((X_1,\ldots,X_{i-1})\) gives
\[
  \E\left\{
    \left\|
      \sum_{j<i}\bm h_0(X_i,X_j)
    \right\|_2^4
    \Bigm|X_i
  \right\}
  \le
  C_d\{(i-1)+(i-1)^2\}
  \E\left(
    \norm{\bm h_0(X_i,X_1)}_2^4
    \mid X_i
  \right).
\]
Taking expectations gives $\E\norm{\bm Y_i}_2^4\le C(i-1)^{-2}$, so $\sum_{i=2}^\infty\E\norm{\bm Y_i}_2^4<\infty$.

Moreover, by Tonelli's theorem,
\[
\begin{aligned}
  \E\left\{
    \sum_{i=2}^{\infty}
    \E(
      \norm{\bm Y_i}_2^4
      \mid\calF_{i-1}
    )
  \right\}
  &=
  \sum_{i=2}^{\infty}
  \E\left\{
    \E(
      \norm{\bm Y_i}_2^4
      \mid\calF_{i-1}
    )
  \right\}=
  \sum_{i=2}^{\infty}
  \E\norm{\bm Y_i}_2^4
  <\infty.
\end{aligned}
\]
Thus $\sum_{i=2}^{\infty}
  \E(
    \norm{\bm Y_i}_2^4
    \mid\calF_{i-1}
  )
  <\infty$ almost surely.

We next establish the covariance-clock bound. Fix
\(\bm a\in\R^d\), let
\(h_{\bm a}(x,y)=\bm a\trans\bm h_0(x,y)\), and define $ q_{\bm a}(x,y)
  =
  \E\{
    h_{\bm a}(Z,x)h_{\bm a}(Z,y)
  \}$ where \(Z\sim P\) is independent. The kernel \(q_{\bm a}\) is
symmetric. It is also canonical, as for \(P\)-a.e. \(y\),
\[
  \E q_{\bm a}(X_1,y)
  =
  \E\left[
    h_{\bm a}(Z,y)
    \E\{h_{\bm a}(Z,X_1)\mid Z\}
  \right]
  =0.
\]
Also, we have
\[
  \E q_{\bm a}(X_1,X_1)
  =
  \E h_{\bm a}(X_1,X_2)^2
  =
  \bm a\trans\bm\Omega\bm a.
\]

Assumption (A2) ensures that both the diagonal and
off-diagonal second moments of \(q_{\bm a}\) are finite. Indeed,
conditional Jensen's inequality gives
\[
  \E q_{\bm a}(X_1,X_1)^2
  \le
  \E h_{\bm a}(Z,X_1)^4
  \le
  \norm{\bm a}_2^4
  \E\norm{\bm h_0(X_1,X_2)}_2^4
  <\infty.
\]
Similarly, applying Jensen's inequality first with respect to \(Z\), we obtain
\[
\begin{aligned}
  \E q_{\bm a}(X_1,X_2)^2
  &\le
  \E\{
    h_{\bm a}(Z,X_1)^2
    h_{\bm a}(Z,X_2)^2
  \} \\
  &=
  \E\left(
    \left[
      \E\bigl\{
        h_{\bm a}(Z,X_1)^2
        \mid Z
      \bigr\}
    \right]^2
  \right) \\
  &\le
  \E h_{\bm a}(Z,X_1)^4
  \le
  \norm{\bm a}_2^4
  \E\norm{\bm h_0(X_1,X_2)}_2^4
  <\infty.
\end{aligned}
\]

Write \(V_{\bm a,n}=\bm a\trans\bracket{\bm M}_n\bm a\) and
\(\mu_{\bm a}=\bm a\trans\bm\Omega\bm a\). Then
\[
\begin{aligned}
  V_{\bm a,n}
  &=
  \sum_{i=2}^n
  \E\left\{
    (\bm a\trans\bm Y_i)^2
    \mid\calF_{i-1}
  \right\} \\
  &=
  \sum_{i=2}^n
  \frac1{(i-1)^2}
  \sum_{1\le j,k<i}q_{\bm a}(X_j,X_k) \\
  &=
  \sum_{j,k=1}^{n-1}
  w_{jk}^{(n)}q_{\bm a}(X_j,X_k).
\end{aligned}
\]
where $w_{jk}^{(n)} =\sum_{r=\max(j,k)}^{n-1}1/{r^2}$. Only the diagonal terms contribute to the expectation, so
\(\E V_{\bm a,n}=\mu_{\bm a}T_n\). Moreover,
\[
  V_{\bm a,n}-\mu_{\bm a}T_n
  =
  \sum_{j=1}^{n-1}
  w_{jj}^{(n)}
  \{q_{\bm a}(X_j,X_j)-\mu_{\bm a}\}
  +
  2\sum_{1\le j<k\le n-1}
  w_{jk}^{(n)}q_{\bm a}(X_j,X_k).
\]
Consider two
distinct unordered pairs \(\{j,k\}\ne\{r,s\}\), with \(j<k\) and
\(r<s\). If the two pairs are disjoint, the independence of $X_i$'s and
\(\E q_{\bm a}(X_j,X_k)=0\) imply $\E\{
    q_{\bm a}(X_j,X_k)
    q_{\bm a}(X_r,X_s)
  \}
  =0$;
If the pairs share exactly one index, for example \(j=r\) and
\(k\ne s\), then
\[
\begin{aligned}
  \E\{
    q_{\bm a}(X_j,X_k)
    q_{\bm a}(X_j,X_s)
  \}=
  \E\left[
    \E\{
      q_{\bm a}(X_j,X_k)
      \mid X_j
    \}
    \E\{
      q_{\bm a}(X_j,X_s)
      \mid X_j
    \}
  \right]
  =0.
\end{aligned}
\]
Thus distinct off-diagonal terms are orthogonal in \(L^2\).

For \(j\ne l\), by canonicity
\[
\begin{aligned}
  \E\left[
    \{q_{\bm a}(X_j,X_j)-\mu_{\bm a}\}
    q_{\bm a}(X_j,X_l)
  \right] =
  \E\left[
    \{q_{\bm a}(X_j,X_j)-\mu_{\bm a}\}
    \E\{
      q_{\bm a}(X_j,X_l)
      \mid X_j
    \}
  \right]
  =0.
\end{aligned}
\]
Consequently,
\[
\begin{aligned}
  \Var(V_{\bm a,n})
  =
  \sum_{j=1}^{n-1}
  \{w_{jj}^{(n)}\}^2
  \Var\{q_{\bm a}(X_j,X_j)\}+
  4\sum_{1\le j<k\le n-1}
  \{w_{jk}^{(n)}\}^2
  \E q_{\bm a}(X_j,X_k)^2.
\end{aligned}
\]
Since \(w_{jj}^{(n)}\le C/j\) and \(w_{jk}^{(n)}\le C/k\) for \(j<k\),
it follows that there exists constant $C_{\bm a}$ such that
\[
  \Var(V_{\bm a,n})
  \le
  C_{\bm a}\sum_{j=1}^{n-1}\frac1{j^2}
  +
  C_{\bm a}\sum_{k=2}^{n-1}\frac{k-1}{k^2}
  \le
  C_{\bm a}T_n.
\]

Let \(n_\ell=\lfloor\exp(\ell^3)\rfloor\). Since
\(T_{n_\ell}\asymp\ell^3\), Chebyshev's inequality gives
\[
  \sum_{\ell=1}^{\infty}
  \Pp\left(
    |V_{\bm a,n_{\ell}}-\mu_{\bm a}T_{n_\ell}|
    >
    T_{n_\ell}^{3/4}
  \right)
  <\infty.
\]
The Borel--Cantelli lemma therefore yields
\(V_{\bm a,n_{\ell}}-\mu_{\bm a}T_{n_\ell}
=O_{a.s.}(T_{n_\ell}^{3/4})\).

For \(n_\ell\le n\le n_{\ell+1}\), monotonicity of \(V_{\bm a,n}\) and
\(T_n\) gives
\[
  V_{\bm a,n}-\mu_{\bm a}T_n
  \le
  V_{\bm a,n_{\ell+1}}-\mu_{\bm a}T_{n_{\ell+1}}
  +
  \mu_{\bm a}(T_{n_{\ell+1}}-T_n)
\]
and
\[
  V_{\bm a,n}-\mu_{\bm a}T_n
  \ge
  V_{\bm a,n_{\ell}}-\mu_{\bm a}T_{n_\ell}
  -
  \mu_{\bm a}(T_n-T_{n_\ell}).
\]
As \(T_{n_{\ell+1}}-T_{n_\ell}
=O(\ell^2)=o(T_{n_\ell}^{3/4})\), it follows that $ V_{\bm a,n}-\mu_{\bm a}T_n
  =
  O_{a.s.}(T_n^{3/4})$ for every \(n\).

Apply this result to the finite collection $\{\bm e_r:1\le r\le d\}
  \cup
  \{\bm e_r+\bm e_s:1\le r<s\le d\}$. Polarization shows that every element of
\(\bracket{\bm M}_n-T_n\bm\Omega\) is
\(O_{a.s.}(T_n^{3/4})\). Since \(d\) is fixed, we have
\[
  \norm{
    \bracket{\bm M}_n-T_n\bm\Omega
  }_{\op}
  =
  O_{a.s.}(T_n^{3/4}),
\]
which also implies
\[
  \norm{
    \bracket{\bm M}_n-T_n\bm\Omega
  }_{\op}
  =
  o_{a.s.}\left(\frac{T_n}{\log T_n}\right),
\]
and division by \(T_n\) proves
\(T_n^{-1}\bracket{\bm M}_n\asto\bm\Omega\) in operator norm.

The scalar optional quadratic-variation result applied to
\(\bm e_r\trans\bm M_n\) and
\((\bm e_r+\bm e_s)\trans\bm M_n\), followed by the same finite
polarization argument, gives $
  \sum_{i=2}^n
  \bm Y_i\bm Y_i\trans /T_n
  \asto\bm\Omega$ in operator norm. Recall that the estimator is
\[
  \widehat{\bm\Omega}_n
  =
  \frac1{T_n}
  \sum_{i=2}^n
  \bm Y_i\bm Y_i\trans
  -
  \frac{
    \bm M_n\bm M_n\trans
  }{
    (n-1)T_n
  }.
\]
For every coordinate \(r\), the scalar ASIP and the Brownian LIL
give
\(M_{n,r}^2=O_{a.s.}\{T_n\log\log T_n\}\). Since \(d\) is fixed,
\(\norm{\bm M_n}_2^2
=O_{a.s.}\{T_n\log\log T_n\}\), and hence
\[
  \frac{
    \norm{\bm M_n\bm M_n\trans}_{\op}
  }{
    (n-1)T_n
  }
  =
  O_{a.s.}\left(
    \frac{\log\log T_n}{n}
  \right)
  \longrightarrow0.
\]
The plug-in consistency now follows from the optional
quadratic-variation convergence.
\end{proof}

\subsection{Strong approximation}

\begin{lemma}[Theorem~1.1 of \citet{zhangASIP}]
\label{lem:zhang-vector-asip}
Let \(\{(\bm\xi_i,\mathcal G_i):i\ge1\}\) be a square-integrable
\(\R^d\)-valued martingale-difference sequence, and denote $ \bm\Sigma_n
  =\sum_{i=1}^n \E( \bm\xi_i\bm\xi_i\trans \mid\mathcal G_{i-1} )
$ and $ V_n=\operatorname{tr}(\bm\Sigma_n).$ Let $f$ be a nondecreasing positive function tending to infinity such that
\(f(t)(\log t)^{\alpha_0}/t\) is nonincreasing for some
\(\alpha_0>50d\), and \(f(t)/t^{\varepsilon_0}\) is nondecreasing
for some \(\varepsilon_0\in(0,1)\). Suppose that \(V_n\to\infty\)
a.s. and
\[
  \sum_{n\ge 1}
    \E\left[
      \norm{\bm\xi_n}_2^2
      I\{\norm{\bm\xi_n}_2^2\ge f(V_n)\}/ f(V_n)
      \mid\mathcal G_{n-1}
    \right] <\infty
  \qquad\text{a.s.},
\]
Moreover, assume that there exists some covariance matrix \(\bm\Gamma\), measurable
with respect to \(\mathcal G_k\) for some \(k\ge0\), such that
\[
  \norm{
    \bm\Sigma_n-\bm\Gamma V_n
  }_{\max}
  =
  o_{a.s.}\{f(V_n)\}.
\]
Then there exists a
sequence \(\{\bm Z_n:n\ge1\}\) of i.i.d. $\mathbb{R}^d$-valued standard normal random vectors, independent of
\(\bm\Gamma\), such that
\[
  \sum_{m\ge1}
  \bm\xi_m I\{V_m\le n\}
  -
  \bm\Gamma^{1/2}
  \sum_{m\le n}\bm Z_m
  =
  O_{a.s.}\left[
    n^{1/2}
    \left\{
      {f(n)}/n
    \right\}^{\frac{1}{50d}}
  \right].
\]
\end{lemma}

\begin{proof}[Proof of Proposition~\ref{prop:vector-asip}]
By the vector setup and Lemma~\ref{lem:vector-clock}, \(\bm Y_i\) is a
square-integrable \(\R^d\)-valued martingale-difference sequence.
Apply Lemma~\ref{lem:zhang-vector-asip} with \(\calG_0\) trivial,
\(\calG_i=\calF_i\) for \(i\ge1\), \(\bm\xi_1=\bm0\), and
\(\bm\xi_i=\bm Y_i\) for \(i\ge2\).
Then \(\bm\Sigma_n=\bracket{\bm M}_n\) and
\(V_n=\operatorname{tr}(\bm\Sigma_n)\). Take
\(c=\operatorname{tr}(\bm\Omega)\) and \(\bm\Gamma={\bm\Omega}/{c}\), then
\[
  \bm\Sigma_n-\bm\Gamma V_n
  =
  \{\bm\Sigma_n-\bm\Omega T_n\}
  -
  \bm\Gamma
  \operatorname{tr}\{\bm\Sigma_n-\bm\Omega T_n\}.
\]
Hence, since \(d\) is fixed,
\[
\begin{aligned}
  \norm{\bm\Sigma_n-\bm\Gamma V_n}_{\max}
  &\le
  \norm{\bm\Sigma_n-\bm\Omega T_n}_{\max}
  +
  \norm{\bm\Gamma}_{\max}
  \left|
    \operatorname{tr}\{\bm\Sigma_n-\bm\Omega T_n\}
  \right| \\
  &\le
  \{1+d\norm{\bm\Gamma}_{\max}\}
  \norm{\bm\Sigma_n-\bm\Omega T_n}_{\op}
  =
  O_{a.s.}(T_n^{3/4}).
\end{aligned}
\]
Taking traces in Lemma~\ref{lem:vector-clock} gives \(V_n/T_n\to c>0\) almost surely, and therefore
\[
  \norm{\bm\Sigma_n-\bm\Gamma V_n}_{\max}
  =
  O_{a.s.}(V_n^{3/4}).
\]
Choose $ f(t)=t^{7/8}$ for all sufficiently large \(t\). Then we obtain
\[
  \norm{\bm\Sigma_n-\bm\Gamma V_n}_{\max}
  =
  o_{a.s.}\{f(V_n)\}.
\]
For example, with $\alpha_0=51d$ and $\varepsilon_0=1/2$, the two monotonicity conditions in Lemma~\ref{lem:zhang-vector-asip} hold for all sufficiently large $t$; extend $f$ positively over the remaining bounded interval. For the conditional Lindeberg condition, \(f(V_i)\ge1\) eventually,
and
\[
  \frac{
    \E\left[
      \norm{\bm Y_i}_2^2
      I\{\norm{\bm Y_i}_2^2\ge f(V_i)\}
      \mid\calF_{i-1}
    \right]
  }{
    f(V_i)
  }
  \le
  \frac{
    \E(
      \norm{\bm Y_i}_2^4
      \mid\calF_{i-1}
    )
  }{
    f(V_i)^2
  }
  \le
  \E(
    \norm{\bm Y_i}_2^4
    \mid\calF_{i-1}
  ).
\]
The series on the right is finite almost surely by
Lemma~\ref{lem:vector-clock}. Thus all conditions of
Lemma~\ref{lem:zhang-vector-asip} are satisfied.

Define
\[
  \nu(m)
  =
  \max\bigl(
    \{1\}\cup\{r\ge2:V_r\le m\}
  \bigr).
\]
There exists a sequence of independent standard normal vectors
\(\{\bm Z_\ell\}\) such that
\[
  \|\bm M_{\nu(m)}
  -
  \bm\Gamma^{1/2}
  \sum_{\ell=1}^m\bm Z_\ell\|_2
  =
  O_{a.s.}(m^{1/2-\kappa_0}),
  \qquad
  \kappa_0=\frac{1-7/8}{50d}=\frac{1}{400d}.
\]
Enlarging the probability space if necessary, construct a standard
\(d\)-dimensional Brownian motion \(\bm B\) satisfying
\(\bm B(m)=\sum_{\ell=1}^m\bm Z_\ell\) for every integer \(m\).
Thus
\[
  \|\bm M_{\nu(m)}-
  \bm\Gamma^{1/2}\bm B(m)\|_2
  =
  O_{a.s.}(m^{1/2-\kappa_0}).
\]

Write
\[
  v_i
  =
  V_i-V_{i-1}
  =
  \E(
    \norm{\bm Y_i}_2^2
    \mid\calF_{i-1}
  ).
\]
Since \(V_{r+1}\) is \(\calF_r\)-measurable,
\(\{\nu(m)\le r\}=\{V_{r+1}>m\}\in\calF_r\); hence \(\nu(m)\)
is a stopping time. Moreover, $ I\{\nu(m-1)<i\le\nu(m)\}
  =
  I\{m-1<V_i\le m\}$ is \(\calF_{i-1}\)-measurable.

Let $ D_m
  =
  \max_{\nu(m-1)\le r\le\nu(m)}
  \norm{
    \bm M_r-\bm M_{\nu(m-1)}
  }_2,$ and $ Q_m = V_{\nu(m)}-V_{\nu(m-1)}.$ If the \(m\)-th block is nonempty, then
\(Q_m\le1+v_{\nu(m-1)+1}\). Conditional Jensen's inequality gives
\(v_i^2\le
\E(\norm{\bm Y_i}_2^4\mid\calF_{i-1})\), and therefore
\[
  \sup_m\E Q_m^2
  \le
  2+
  2\sum_{i=2}^{\infty}
  \E\norm{\bm Y_i}_2^4
  <\infty.
\]
Since
\(I\{\nu(m-1)<i\le\nu(m)\}\) is
\(\calF_{i-1}\)-measurable, write \(Y_{i,r}\) for the \(r\)-th coordinate of
\(\bm Y_i\). Then the stopped increments
\[
  Y_{i,r}I\{\nu(m-1)<i\le\nu(m)\},
  \qquad i\ge2,
\]
form a scalar martingale-difference sequence for each
\(r=1,\ldots,d\).

Using \(\norm{\bm x}_2^4\le d\sum_{r=1}^d|x_r|^4\), and applying Doob's
\(L^4\) maximal inequality followed by the scalar martingale
Rosenthal inequality
\citep[Theorem~2.12]{HallHeyde1980} coordinatewise, we obtain
\[
\begin{aligned}
  \E D_m^4
  \le C_d\sum_{r=1}^d\E\Bigg[
    &
    \left\{
      \sum_{i=\nu(m-1)+1}^{\nu(m)}
      \E(Y_{i,r}^2\mid\calF_{i-1})
    \right\}^2 \\
    &+
    \sum_{i=\nu(m-1)+1}^{\nu(m)}
    \E(|Y_{i,r}|^4\mid\calF_{i-1})
  \Bigg].
\end{aligned}
\]
Since all conditional second moments are nonnegative,
\[
\begin{aligned}
  \sum_{r=1}^d
  \left\{
    \sum_{i=\nu(m-1)+1}^{\nu(m)}
    \E(Y_{i,r}^2\mid\calF_{i-1})
  \right\}^2
  &\le
  \left\{
    \sum_{i=\nu(m-1)+1}^{\nu(m)}
    \E(\norm{\bm Y_i}_2^2\mid\calF_{i-1})
  \right\}^2 \\
  &=Q_m^2.
\end{aligned}
\]
Moreover,
\[
  \sum_{r=1}^d
  \sum_{i=\nu(m-1)+1}^{\nu(m)}
  \E(|Y_{i,r}|^4\mid\calF_{i-1})
  \le
  \sum_{i=\nu(m-1)+1}^{\nu(m)}
  \E(\norm{\bm Y_i}_2^4\mid\calF_{i-1}).
\]
Consequently,
\[
  \E D_m^4
  \le
  C_d\E\left[
    Q_m^2+
    \sum_{i=\nu(m-1)+1}^{\nu(m)}
    \E(\norm{\bm Y_i}_2^4\mid\calF_{i-1})
  \right].
\]
The preceding bound on \(Q_m\) and
Lemma~\ref{lem:vector-clock} therefore imply
\[
  \sup_{m\ge1}\E D_m^4<\infty.
\]
Hence, for every \(\eta>0\),
\[
  \sum_{m=1}^\infty
  \Pp\{D_m>m^{1/4+\eta}\}
  \le
  C\sum_{m=1}^\infty m^{-1-4\eta}
  <\infty,
\]
and the Borel--Cantelli lemma gives $D_m=O_{a.s.}(m^{1/4+\eta}).$

For \(m_n=\lceil V_n\rceil\), we have
\(m_n-1<V_n\le m_n\), and hence
\(\nu(m_n-1)<n\le\nu(m_n)\). Therefore,
\[
\begin{aligned}
  \norm{
    \bm M_n-\bm M_{\nu(m_n)}
  }_2
  &\le
  \norm{
    \bm M_n-\bm M_{\nu(m_n-1)}
  }_2
  +
  \norm{
    \bm M_{\nu(m_n)}
    -
    \bm M_{\nu(m_n-1)}
  }_2 \\
  &\le 2D_{m_n}
  =
  O_{a.s.}(V_n^{1/4+\delta})
\end{aligned}
\]
for every \(\delta>0\), where we used \(m_n/V_n\to1\) almost
surely.

Moreover, the Brownian unit-interval oscillations satisfy
\[
  \sup_{m-1\le t\le m}
  \norm{\bm B(t)-\bm B(m)}_2
  =
  O_{a.s.}(\sqrt{\log m}).
\]
Combining this bound with the approximation at \(\nu(m_n)\) gives
\[
\begin{aligned}
  \norm{
    \bm M_n-\bm\Gamma^{1/2}\bm B(V_n)
  }_2
  =
  O_{a.s.}\left(
    V_n^{1/4+\delta}
    +
    V_n^{1/2-\kappa_0}
    +
    \sqrt{\log V_n}
  \right).
\end{aligned}
\]
Choosing \(0<\delta<1/4\), we obtain
\[
  \norm{
    \bm M_n-\bm\Gamma^{1/2}\bm B(V_n)
  }_2
  =
  O_{a.s.}(V_n^{1/2-\kappa_1}),
  \qquad
  \kappa_1
  =
  \min\left\{
    \kappa_0,\frac14-\delta
  \right\}>0.
\]

Finally, Lemma~\ref{lem:vector-clock} gives
\[
  |V_n-cT_n|
  =
  O_{a.s.}(T_n^{3/4}).
\]
A standard Brownian modulus estimate on dyadic time blocks implies
that, for every \(\varepsilon>0\),
\[
  \norm{
    \bm B(V_n)-\bm B(cT_n)
  }_2
  =
  O_{a.s.}(T_n^{3/8+\varepsilon}).
\]
Since \(V_n\asymp T_n\) almost surely, choosing
\(0<\varepsilon<1/8\) gives
\[
  \norm{
    \bm M_n-\bm\Gamma^{1/2}\bm B(cT_n)
  }_2
  =
  O_{a.s.}(T_n^{1/2-\kappa}),
\]
where
\[
  \kappa
  =
  \min\left\{
    \kappa_1,\frac18-\varepsilon
  \right\}>0.
\]

Define
\(\bm W(t)=c^{-1/2}\bm B(ct)\). By Brownian scaling,
\(\bm W\) is a standard \(d\)-dimensional Brownian motion, and since
\(\bm\Gamma=\bm\Omega/c\),
\[
  \bm\Gamma^{1/2}\bm B(cT_n)
  =
  \bm\Omega^{1/2}\bm W(T_n).
\]
This completes the proof.
\end{proof}

\subsection{Radial Brownian boundaries}

\begin{lemma}[Vector radial Brownian boundaries]
\label{lem:vector-brownian-boundaries}
For either $b_{d,\alpha,t_0}
    \in
    \left\{
        b^{\mathrm{GM}}_{d,\alpha,t_0},
        b^{\mathrm{LIL}}_{d,\alpha,t_0}
    \right\},$
the bound satisfies
\begin{equation}
    \Pp\left\{
        \exists t\ge t_0:
        \norm{\bm W_d(t)}_2
        \ge
        b_{d,\alpha,t_0}(t)
    \right\}
    \le \alpha.
    \label{eq:vector-boundary-crossing}
\end{equation}
Equally, both choices are radial
$(1-\alpha)$ Brownian motion boundaries in the sense of
Definition~\ref{def:vector-boundary}. Moreover, for every fixed $t_0>0$,
\begin{equation}
    \limsup_{t\to\infty}
    \frac{\norm{\bm W_d(t)}_2}
         {b_{d,\alpha,t_0}(t)}
    <1
    \qquad\text{a.s.}
    \label{eq:vector-boundary-tail}
\end{equation}
\end{lemma}

\begin{proof}
For the Gaussian-mixture boundary, let
$\mu=\mathcal N_d(0,t_0^{-1}I_d)$ and consider
\[
\begin{aligned}
L_t&=\int_{\mathbb R^d}
\exp\left\{\lambda^\top\bm W_d(t)
-\frac{t\|\lambda\|_2^2}{2}\right\}\,\mu(d\lambda) \\
&=\left(\frac{t_0}{t+t_0}\right)^{d/2}
\exp\left\{\frac{\|\bm W_d(t)\|_2^2}{2(t+t_0)}\right\},
\qquad t\ge0.
\end{aligned}
\]
For each fixed $\lambda$, the integrand is a mean-one
exponential martingale. Since $\mu$ is a probability measure,
conditional Tonelli's theorem shows that $\{L_t\}_{t\ge0}$
is a nonnegative martingale with $L_0=1$.
For $t\ge t_0$, a crossing of
$b^{\mathrm{GM}}_{d,\alpha,t_0}$ implies $L_t\ge1/\alpha$.
Ville's inequality therefore gives
\[
\Pp\left\{\exists t\ge t_0:
\|\bm W_d(t)\|_2\ge b^{\mathrm{GM}}_{d,\alpha,t_0}(t)
\right\}
\le \Pp\left\{\sup_{t\ge0}L_t\ge1/\alpha\right\}
\le\alpha.
\]
This proves \eqref{eq:vector-boundary-crossing} for the
Gaussian-mixture boundary.

For the LIL-type boundary, for every
$v\in\mathbb S^{d-1}$ and $\lambda\in\R$, $\exp\left\{
        \lambda v\trans\bm W_d(t)
        -{\lambda^2t}/{2}
    \right\}$ is a nonnegative martingale for $t\ge 0$. Hence
$\{(\bm W_d(t),tI_d)\}_{t\ge0}$ is a vector sub-$\psi_N$
process with $\psi_N(\lambda)=\lambda^2/2$.

Apply \citet[Theorem~4.1]{WhitehouseWuRamdas2025} with
regularization parameter $\rho=t_0$, confidence level
$\delta=\alpha$, and $h(x)=\zeta(q)(x+1)^q$. Since $\sum_{k=0}^{\infty}1/{h(k)}=1$, the theorem applies.

For $t\ge t_0$, $tI_d\vee t_0I_d=tI_d,
    \gamma_{\min}(tI_d)
    =
    \gamma_{\max}(tI_d)
    =t,
    \kappa(tI_d)=1,$
and $(\psi_N^*)^{-1}(x)=\sqrt{2x}.$
Therefore, Theorem~4.1 yields that
\[
    \Pp\left\{
        \exists t\ge t_0:
        \norm{\bm W_d(t)}_2
        \ge
        \frac{\sqrt{2at\,L_{t_0}(tI_d)}}{1-\varepsilon}
    \right\}
    \le\alpha,
\]
where $L_{t_0}(tI_d)
    =
    \log\zeta(q)
    +
    q\log\left\{
        1+\log_a\left({t}/{t_0}\right)
    \right\} +
    \log\{{1}/{\alpha(1-\beta^{-1})}\}
    +
    \log\left\{
        \beta
        N_{d-1}\left({\varepsilon}/{\beta}\right)
    \right\}.$ Using $N_{d-1}(u)
    \le
    \left(3/u\right)^d,$ we obtain
\[
\begin{aligned}
    L_{t_0}(tI_d)
    \le{}&
    q\log\left\{
        1+\log_a\left(\frac{t}{t_0}\right)
    \right\}
    +\log\zeta(q) \\
    &+
    \log\frac{1}{\alpha(1-\beta^{-1})}
    +\log\beta
    +d\log\left(\frac{3\beta}{\varepsilon}\right).
\end{aligned}
\]
Enlarging the boundary therefore gives
\eqref{eq:vector-boundary-crossing} with
$b_{d,\alpha,t_0}=b^{\mathrm{LIL}}_{d,\alpha,t_0}$.

From the \eqref{eq:vector-gm-boundary} and \eqref{eq:vector-lil-boundary}, we know that $b_{d,\alpha,t_0}(t_0r)
    =
    \sqrt{t_0}\,
    b_{d,\alpha,1}(r)$ for $r\ge1$. Moreover, for the Gaussian-mixture boundary, since $t\ge t_0$, $b^{\mathrm{GM}}_{d,\alpha,t_0}(t)
    \ge
    \sqrt{
        d\log2+2\log(1/\alpha)
    }\sqrt t.$ For the LIL-type boundary, the expression inside the square brackets
in \eqref{eq:vector-lil-boundary} is bounded below by $\log\zeta(q)
    +
    \log\{{1}/{\alpha(1-\beta^{-1})}\}
    +
    \log\beta
    +
    d\log\left({3\beta}/{\varepsilon}\right)
    >0.$ Thus there exists a constant $c_0>0$, independent of $t_0$, such that
\begin{equation*}
    b_{d,\alpha,t_0}(t)
    \ge c_0\sqrt{t},
    \qquad t\ge t_0>0.
\end{equation*}

Finally, for fixed $t_0$, ${b^{\mathrm{GM}}_{d,\alpha,t_0}(t)}/{\sqrt{2t\log\log t}}\rightarrow\infty,$
whereas ${
        b^{\mathrm{LIL}}_{d,\alpha,t_0}(t)
    }/{
        \sqrt{2t\log\log t}
    }
    \rightarrow
    {\sqrt{aq}}/{(1-\varepsilon)}
    >1$.
The fixed-dimensional Brownian law of the iterated logarithm gives
\[
    \limsup_{t\to\infty}
    \frac{\norm{\bm W_d(t)}_2}
         {\sqrt{2t\log\log t}}
    =1
    \qquad\text{a.s.}
\]
Hence \eqref{eq:vector-boundary-tail} follows for both boundaries.
\end{proof}

\subsection{Time-uniform coverage}

\begin{proof}[Proof of Theorem~\ref{thm:vector-coverage}]
By Proposition~\ref{prop:vector-asip}, on an enriched probability space
supporting a standard \(d\)-dimensional Brownian motion \(\bm W_d\), $\bm\Omega^{-1/2}\bm M_n
  =
  \bm W_d(T_n)+\bm e_n$ where $\norm{\bm e_n}_2
  =
  o_{a.s.}(\sqrt{T_n})$.
By Lemma~\ref{lem:vector-clock},
\(\widehat{\bm\Omega}_n\asto\bm\Omega\).
Since \(\bm\Omega\) is positive definite,
\(\widehat{\bm\Omega}_n\) is almost surely positive definite for all
sufficiently large \(n\).

For every \(n\), define
\[
  \bm A_n
  =
  \begin{cases}
    \bm\Omega^{1/2}
    \widehat{\bm\Omega}_n^{-1}
    \bm\Omega^{1/2},
    &\widehat{\bm\Omega}_n\succ0,\\[4pt]
    \bm0_{d\times d},
    &\widehat{\bm\Omega}_n\text{ is singular}.
  \end{cases}
\]
Then $\norm{\bm A_n-\bm I_d}_{\rm op}\asto0.$

By Lemma~\ref{lem:vector-brownian-boundaries}, for some \(c_0>0\), $b_{d,\alpha,T_m}(T_n)\ge c_0\sqrt{T_n}$ for all $n\ge m$. Define $$\eta_m
  =
  1-
  \left\{
    1+
    \sup_{n\ge m}
    \norm{\bm A_n-\bm I_d}_{\rm op}
  \right\}^{-1/2} +\frac{1}{c_0}
  \sup_{n\ge m}\frac{\norm{\bm e_n}_2}{\sqrt{T_n}}.$$ Then $\eta_m\asto0.$

Suppose that, for some \(n\ge m\),
\(\bm\theta\notin C_{m,n}^{\rm ell}(\alpha)\).
By the convention in \eqref{eq:vector-ellipsoid},
this can occur only if \(\widehat{\bm\Omega}_n\succ0\), in which case
\[
  \bm M_n\trans
  \widehat{\bm\Omega}_n^{-1}
  \bm M_n
  >
  \left\{
    b_{d,\alpha,T_m}(T_n)
  \right\}^2.
\]
Using Proposition \ref{prop:vector-asip},
\[
\begin{aligned}
  b_{d,\alpha,T_m}(T_n)
  &<
  \left[
    \left\{
      \bm W_d(T_n)+\bm e_n
    \right\}\trans
    \bm A_n
    \left\{
      \bm W_d(T_n)+\bm e_n
    \right\}
  \right]^{1/2} \\
  &\le
  \sqrt{
    1+\norm{\bm A_n-\bm I_d}_{\rm op}
  }\,
  \left\{
    \norm{\bm W_d(T_n)}_2
    +
    \norm{\bm e_n}_2
  \right\}.
\end{aligned}
\]
Therefore
\[
\begin{aligned}
  \norm{\bm W_d(T_n)}_2
  &>
  \left\{
    1+\norm{\bm A_n-\bm I_d}_{\rm op}
  \right\}^{-1/2}
  b_{d,\alpha,T_m}(T_n)
  -
  \norm{\bm e_n}_2 \\
  &\ge
  (1-\eta_m)
  b_{d,\alpha,T_m}(T_n).
\end{aligned}
\]

Define the feasible crossing events $E_m
  =
  \left\{
    \exists n\ge m:
    \bm\theta\notin C_{m,n}^{\rm ell}(\alpha)
  \right\}$ and ideal crossing events $A_m
  =
  \left\{
    \exists t\ge T_m:
    \norm{\bm W_d(t)}_2
    \ge
    b_{d,\alpha,T_m}(t)
  \right\}.$
By Lemma~\ref{lem:vector-brownian-boundaries}, $\Pp(A_m)\le\alpha.$

Fix \(\delta\in(0,1)\). On \(\{\eta_m\le\delta\}\), occurrence of
\(E_m\) implies that, for some \(n\ge m\),
\[
  \norm{\bm W_d(T_n)}_2
  >
  (1-\delta)
  b_{d,\alpha,T_m}(T_n).
\]
Hence
\[
\begin{aligned}
  \Pp(E_m\setminus A_m)
  \le
  \Pp(\eta_m>\delta)+
  \Pp\left\{
    A_m^c,\
    \exists t\ge T_m:
    \norm{\bm W_d(t)}_2
    >
    (1-\delta)b_{d,\alpha,T_m}(t)
  \right\}.
\end{aligned}
\]
Since \(\eta_m\asto0\), $\Pp(\eta_m>\delta)\to0.$

By $b_{d,\alpha,t_0}(t_0r)=
    \sqrt{t_0}\,
    b_{d,\alpha,1}(r)$ and Brownian scaling, the second
probability equals
\[
  \Pp\left\{
    \begin{gathered}
    \norm{\bm W_d(r)}_2 < b_{d,\alpha,1}(r)
    \text{ for all }r\ge1,\\
    \exists r\ge1:\quad
    \norm{\bm W_d(r)}_2 > (1-\delta)b_{d,\alpha,1}(r)
    \end{gathered}
  \right\}.
\]
From Lemma \ref{lem:vector-brownian-boundaries},  \eqref{eq:vector-boundary-tail} gives that
\[
  \limsup_{r\to\infty}
  \frac{\norm{\bm W_d(r)}_2}
       {b_{d,\alpha,1}(r)}
  <1
  \qquad\text{a.s.}
\]
Thus, on the event $\left\{
    \norm{\bm W_d(r)}_2
    <
    b_{d,\alpha,1}(r)
    \text{ for all }r\ge1
  \right\},$ there exist \(\varepsilon>0\) and \(R<\infty\) such that
\[
  \frac{\norm{\bm W_d(r)}_2}
       {b_{d,\alpha,1}(r)}
  \le
  1-\varepsilon,
  \qquad r\ge R.
\]
On the compact interval \([1,R]\), continuity of
\(\bm W_d\) and \(b_{d,\alpha,1}\) gives
\[
  \sup_{1\le r\le R}
  \frac{\norm{\bm W_d(r)}_2}
       {b_{d,\alpha,1}(r)}
  <1.
\]
Consequently,
\[
  \sup_{r\ge1}
  \frac{\norm{\bm W_d(r)}_2}
       {b_{d,\alpha,1}(r)}
  <1
\]
almost surely on this event. It follows that as $\delta\downarrow0$,
$$\Pp\left\{
  \norm{\bm W_d(r)}_2
  <
  b_{d,\alpha,1}(r)
  \text{ for all }r\ge1,\
  \exists r\ge1:
  \norm{\bm W_d(r)}_2
  >
  (1-\delta)b_{d,\alpha,1}(r)
\right\}
\rightarrow0$$

We therefore obtain, for every fixed \(\delta\in(0,1)\),
\[
\begin{aligned}
  \limsup_{m\to\infty}\Pp(E_m)
  &\le
  \alpha
  +
  \Pp\left\{
    \begin{gathered}
    \norm{\bm W_d(r)}_2 < b_{d,\alpha,1}(r)
    \text{ for all }r\ge1,\\
    \exists r\ge1:\quad
    \norm{\bm W_d(r)}_2 > (1-\delta)b_{d,\alpha,1}(r)
    \end{gathered}
  \right\}.
\end{aligned}
\]
Letting \(\delta\downarrow0\) yields
\[
  \limsup_{m\to\infty}
  \Pp\left\{
    \exists n\ge m:
    \bm\theta\notin C_{m,n}^{\rm ell}(\alpha)
  \right\}
  \le\alpha.
\]
Equivalently,
\[
  \liminf_{m\to\infty}
  \Pp\left\{
    \forall n\ge m:
    \bm\theta\in C_{m,n}^{\rm ell}(\alpha)
  \right\}
  \ge1-\alpha.
\]
\end{proof}

\subsection{Multivariate asymptotic confidence sequences}
\begin{proof}[Proof of Remark \ref{rem:vector-asympcs}]
Fix \(m\). By Proposition~\ref{prop:vector-asip},
\[
  \bm\Omega^{-1/2}\bm M_n
  =
  \bm W_d(T_n)+\bm e_n,
  \qquad
  \norm{\bm e_n}_2
  =
  o_{a.s.}(\sqrt{T_n}).
\]
Hence
\[
  \sup_{n\ge m}
  \frac{\norm{\bm e_n}_2}{\sqrt{T_n}}
  <\infty
  \qquad\text{a.s.}
\]
Choose \(\delta\in(0,\alpha)\) and a deterministic
\(\varepsilon>0\) such that
\[
  \Pp\left\{
    \sup_{n\ge m}
    \frac{\norm{\bm e_n}_2}{\sqrt{T_n}}
    >
    \varepsilon
  \right\}
  \le\delta,
\]
and put \(a=\alpha-\delta\). Define
\[
  C_{m,n}^{*,\rm ell}(\alpha)
  =
  \left\{
    \bm\vartheta\in\mathbb R^d:
    \left\|
      \bm\Omega^{-1/2}
      (n-1)
      (\widetilde{\bm U}_n-\bm\vartheta)
    \right\|_2
    \le
    b_{d,a,T_m}(T_n)
    +
    \varepsilon\sqrt{T_n}
  \right\}.
\]
By Lemma~\ref{lem:vector-brownian-boundaries},
\[
  \Pp\left\{
    \norm{\bm W_d(t)}_2
    <
    b_{d,a,T_m}(t)
    \text{ for all }t\ge T_m
  \right\}
  \ge1-a.
\]
On this event and
\[
  \left\{
    \sup_{n\ge m}
    \frac{\norm{\bm e_n}_2}{\sqrt{T_n}}
    \le\varepsilon
  \right\},
\]
we have, for every \(n\ge m\),
\[
\begin{aligned}
  \norm{\bm\Omega^{-1/2}\bm M_n}_2
  &\le
  \norm{\bm W_d(T_n)}_2
  +
  \norm{\bm e_n}_2 \\
  &<
  b_{d,a,T_m}(T_n)
  +
  \varepsilon\sqrt{T_n}.
\end{aligned}
\]
Since $\bm M_n
  =
  (n-1)(\widetilde{\bm U}_n-\bm\theta),$ a union bound gives
\[
  \Pp\left\{
    \bm\theta\in
    C_{m,n}^{*,\rm ell}(\alpha)
    \text{ for all }n\ge m
  \right\}
  \ge
  1-a-\delta
  =
  1-\alpha.
\]
Thus
\(\{C_{m,n}^{*,\rm ell}(\alpha):n\ge m\}\)
is a nonasymptotic \((1-\alpha)\)-confidence sequence.

It remains to compare this sequence with
\(C_{m,n}^{\rm ell}(\alpha)\). For fixed \(m\), the explicit forms of
either radial boundary give
\[
  \frac{
    b_{d,a,T_m}(T_n)
    +
    \varepsilon\sqrt{T_n}
  }{
    b_{d,\alpha,T_m}(T_n)
  }
  \longrightarrow1.
\]
Indeed, as ${
    b_{d,a,T_m}(T_n)
  }/{
    b_{d,\alpha,T_m}(T_n)
  }
  \rightarrow1,$ ${\sqrt{T_n}}/
       {b_{d,\alpha,T_m}(T_n)}\rightarrow0$
and $\widehat{\bm\Omega}_n\asto\bm\Omega$ there exists a random sequence
\(\rho_n\asto0\) such that, almost surely for all sufficiently large
\(n\),
\[
  (1-\rho_n)
  \left\{
    C_{m,n}^{*,\rm ell}(\alpha)
    -
    \widetilde{\bm U}_n
  \right\}
  \subseteq
  C_{m,n}^{\rm ell}(\alpha)
  -
  \widetilde{\bm U}_n
  \subseteq
  (1+\rho_n)
  \left\{
    C_{m,n}^{*,\rm ell}(\alpha)
    -
    \widetilde{\bm U}_n
  \right\}.
\]
Therefore,
\[
\begin{aligned}
  \frac{
    \nu\left(
      C_{m,n}^{\rm ell}(\alpha)
      \mathbin{\triangle}
      C_{m,n}^{*,\rm ell}(\alpha)
    \right)
  }{
    \nu\left(
      C_{m,n}^{*,\rm ell}(\alpha)
    \right)
  }
  &\le
  (1+\rho_n)^d-(1-\rho_n)^d
  \asto0.
\end{aligned}
\]
Hence the ellipsoids satisfy the multivariate AsympCS criterion of
\citet[Definition~B.6]{WaudbySmithArbourSinhaKennedyRamdas2024}.
\end{proof}


\begin{thebibliography}{36}
\providecommand{\natexlab}[1]{#1}
\providecommand{\url}[1]{\texttt{#1}}
\expandafter\ifx\csname urlstyle\endcsname\relax
  \providecommand{\doi}[1]{doi: #1}\else
  \providecommand{\doi}{doi: \begingroup \urlstyle{rm}\Url}\fi

\bibitem[Arcones and Gin{\'e}(1993)]{ArconesGine1993}
Miguel~A. Arcones and Evarist Gin{\'e}.
\newblock Limit theorems for {$U$}-processes.
\newblock \emph{The Annals of Probability}, 21\penalty0 (3):\penalty0
  1494--1542, 1993.

\bibitem[Arcones and Gin{\'e}(1995)]{ArconesGine1995}
Miguel~A. Arcones and Evarist Gin{\'e}.
\newblock On the law of the iterated logarithm for canonical {$U$}-statistics
  and processes.
\newblock \emph{Stochastic Processes and their Applications}, 58\penalty0
  (2):\penalty0 217--245, 1995.

\bibitem[Bickel and Levina(2008)]{BickelLevina2008}
Peter~J. Bickel and Elizaveta Levina.
\newblock Covariance regularization by thresholding.
\newblock \emph{The Annals of Statistics}, 36\penalty0 (6):\penalty0
  2577--2604, 2008.

\bibitem[Cai et~al.(2026)Cai, Hu, and Li]{ati2026}
Leheng Cai, Qirui Hu, and Weijia Li.
\newblock Asymptotic anytime-valid inference for {U}-statistics, 2026.
\newblock URL \url{https://arxiv.org/abs/2605.14692}.

\bibitem[Cattaneo et~al.(2025)Cattaneo, Masini, and
  Underwood]{CattaneoMasiniUnderwood2025}
Matias~D. Cattaneo, Ricardo~P. Masini, and William~G. Underwood.
\newblock {Yurinskii}'s coupling for martingales.
\newblock \emph{The Annals of Statistics}, 53\penalty0 (5):\penalty0
  2179--2203, 2025.

\bibitem[Chatterjee and Ramdas(2026)]{ChatterjeeRamdas2025mMMD}
Anirban Chatterjee and Aaditya Ramdas.
\newblock A martingale kernel two-sample test.
\newblock In \emph{Proceedings of the 37th International Conference on
  Algorithmic Learning Theory}, volume 313 of \emph{Proceedings of Machine
  Learning Research}, pages 1--44, 2026.

\bibitem[de~la Pe{\~n}a et~al.(2007)de~la Pe{\~n}a, Klass, and
  Lai]{de_la_Pe_a_2007}
Victor~H. de~la Pe{\~n}a, Michael~J. Klass, and Tze~Leung Lai.
\newblock Pseudo-maximization and self-normalized processes.
\newblock \emph{Probability Surveys}, 4:\penalty0 172--192, 2007.

\bibitem[Dehling et~al.(1984)Dehling, Denker, and
  Philipp]{DehlingDenkerPhilipp1984}
Herold Dehling, Manfred Denker, and Walter Philipp.
\newblock Invariance principles for von {Mises} and {U}-statistics.
\newblock \emph{Zeitschrift f{\"u}r Wahrscheinlichkeitstheorie und Verwandte
  Gebiete}, 67:\penalty0 139--167, 1984.

\bibitem[Duan et~al.(2022)Duan, Ramdas, and Wasserman]{DuanRamdasWasserman2022}
Boyan Duan, Aaditya Ramdas, and Larry Wasserman.
\newblock Interactive rank testing by betting.
\newblock In \emph{Proceedings of the First Conference on Causal Learning and
  Reasoning}, volume 177 of \emph{Proceedings of Machine Learning Research},
  pages 201--235. PMLR, 2022.

\bibitem[Gnettner and Kirch(2025)]{GnettnerKirch2025}
Felix Gnettner and Claudia Kirch.
\newblock A new and flexible class of sharp asymptotic time-uniform confidence
  sequences.
\newblock \emph{Statistics \& Probability Letters}, 226:\penalty0 110462, 2025.

\bibitem[Gretton et~al.(2005)Gretton, Bousquet, Smola, and
  Sch{\"o}lkopf]{GrettonEtAl2005}
Arthur Gretton, Olivier Bousquet, Alex Smola, and Bernhard Sch{\"o}lkopf.
\newblock Measuring statistical dependence with {Hilbert--Schmidt} norms.
\newblock In \emph{Algorithmic Learning Theory}, volume 3734 of \emph{Lecture
  Notes in Computer Science}, pages 63--77. Springer, 2005.

\bibitem[Gretton et~al.(2012)Gretton, Borgwardt, Rasch, Sch{\"o}lkopf, and
  Smola]{GrettonEtAl2012}
Arthur Gretton, Karsten~M. Borgwardt, Malte~J. Rasch, Bernhard Sch{\"o}lkopf,
  and Alexander Smola.
\newblock A kernel two-sample test.
\newblock \emph{Journal of Machine Learning Research}, 13:\penalty0 723--773,
  2012.

\bibitem[Gr{\"u}nwald et~al.(2024)Gr{\"u}nwald, Henzi, and
  Lardy]{GrunwaldHenziLardy2024}
Peter Gr{\"u}nwald, Alexander Henzi, and Tyron Lardy.
\newblock Anytime-valid tests of conditional independence under {Model-X}.
\newblock \emph{Journal of the American Statistical Association}, 119\penalty0
  (546):\penalty0 1554--1565, 2024.

\bibitem[Hall and Heyde(1980)]{HallHeyde1980}
Peter Hall and Christopher~C. Heyde.
\newblock \emph{Martingale Limit Theory and Its Application}.
\newblock Academic Press, 1980.

\bibitem[Hoeffding(1948)]{Hoeffding1948}
Wassily Hoeffding.
\newblock A class of statistics with asymptotically normal distribution.
\newblock \emph{The Annals of Mathematical Statistics}, 19\penalty0
  (3):\penalty0 293--325, 1948.

\bibitem[Howard et~al.(2020)Howard, Ramdas, McAuliffe, and
  Sekhon]{HowardRamdasMcAuliffeSekhon2020}
Steven~R. Howard, Aaditya Ramdas, Jon McAuliffe, and Jasjeet Sekhon.
\newblock Time-uniform {Chernoff} bounds via nonnegative supermartingales.
\newblock \emph{Probability Surveys}, 17:\penalty0 257--317, 2020.

\bibitem[Howard et~al.(2021)Howard, Ramdas, McAuliffe, and
  Sekhon]{HowardRamdasMcAuliffeSekhon2021}
Steven~R. Howard, Aaditya Ramdas, Jon McAuliffe, and Jasjeet Sekhon.
\newblock Time-uniform, nonparametric, nonasymptotic confidence sequences.
\newblock \emph{The Annals of Statistics}, 49\penalty0 (2):\penalty0
  1055--1080, 2021.

\bibitem[Kanagawa and Yoshihara(1994)]{KanagawaYoshihara1994}
S.~Kanagawa and K.~Yoshihara.
\newblock The almost sure invariance principles of degenerate {U}-statistics of
  degree two for stationary random variables.
\newblock \emph{Stochastic Processes and their Applications}, 49\penalty0
  (2):\penalty0 347--356, 1994.
\newblock ISSN 0304-4149.

\bibitem[Kim and Ramdas(2024)]{KimRamdas2024}
Ilmun Kim and Aaditya Ramdas.
\newblock Dimension-agnostic inference using cross {$U$}-statistics.
\newblock \emph{Bernoulli}, 30\penalty0 (1):\penalty0 683--711, 2024.

\bibitem[Liu et~al.(2026)Liu, Feng, Peng, and Wang]{LiuEtAl2026BootAVT}
\begingroup
Jixuan Liu, Long Feng, Liuhua Peng, and Zhaojun Wang.
\newblock {Boot-AVT}: Asymptotically anytime-valid testing with degenerate
  {U}-statistics via bootstrap, 2026.
\newblock URL \url{https://arxiv.org/abs/2610.28931}.
\newblock arXiv:2610.28931v1.
\par\endgroup

\bibitem[Martinez-Taboada and Ramdas(2025)]{MartinezTaboadaRamdas2025}
Diego Martinez-Taboada and Aaditya Ramdas.
\newblock Sequential kernelized {Stein} discrepancy.
\newblock In \emph{Proceedings of the 28th International Conference on
  Artificial Intelligence and Statistics}, volume 258 of \emph{Proceedings of
  Machine Learning Research}, pages 1288--1296. PMLR, 2025.

\bibitem[Neuhaus(1977)]{Neuhaus1977}
Georg Neuhaus.
\newblock Functional limit theorems for {$U$}-statistics in the degenerate
  case.
\newblock \emph{Journal of Multivariate Analysis}, 7\penalty0 (3):\penalty0
  424--439, 1977.

\bibitem[Podkopaev et~al.(2023)Podkopaev, Bl{\"o}baum, Kasiviswanathan, and
  Ramdas]{PodkopaevEtAl2023}
Aleksandr Podkopaev, Patrick Bl{\"o}baum, Shiva Kasiviswanathan, and Aaditya
  Ramdas.
\newblock Sequential kernelized independence testing.
\newblock In \emph{Proceedings of the 40th International Conference on Machine
  Learning}, volume 202 of \emph{Proceedings of Machine Learning Research},
  pages 27957--27993. PMLR, 2023.

\bibitem[Ramdas et~al.(2023)Ramdas, Gr{\"u}nwald, Vovk, and
  Shafer]{RamdasGrunwaldVovkShafer2023}
Aaditya Ramdas, Peter Gr{\"u}nwald, Vladimir Vovk, and Glenn Shafer.
\newblock Game-theoretic statistics and safe anytime-valid inference.
\newblock \emph{Statistical Science}, 38\penalty0 (4):\penalty0 576--601, 2023.

\bibitem[Robbins(1970)]{Robbins1970a}
Herbert Robbins.
\newblock Statistical methods related to the law of the iterated logarithm.
\newblock \emph{The Annals of Mathematical Statistics}, 41\penalty0
  (5):\penalty0 1397--1409, 1970.

\bibitem[Robbins and Siegmund(1970)]{Robbins1970b}
Herbert Robbins and David Siegmund.
\newblock Boundary crossing probabilities for the {Wiener} process and sample
  sums.
\newblock \emph{The Annals of Mathematical Statistics}, 41\penalty0
  (5):\penalty0 1410--1429, 1970.

\bibitem[Sejdinovic et~al.(2013)Sejdinovic, Sriperumbudur, Gretton, and
  Fukumizu]{SejdinovicEtAl2013}
Dino Sejdinovic, Bharath Sriperumbudur, Arthur Gretton, and Kenji Fukumizu.
\newblock Equivalence of distance-based and {RKHS}-based statistics in
  hypothesis testing.
\newblock \emph{The Annals of Statistics}, 41\penalty0 (5):\penalty0
  2263--2291, 2013.

\bibitem[Serfling(1980)]{Serfling1980}
Robert~J. Serfling.
\newblock \emph{Approximation Theorems of Mathematical Statistics}.
\newblock John Wiley \& Sons, 1980.

\bibitem[Shekhar and Ramdas(2024)]{ShekharRamdas2024}
Shubhanshu Shekhar and Aaditya Ramdas.
\newblock Nonparametric two-sample testing by betting.
\newblock \emph{IEEE Transactions on Information Theory}, 70\penalty0
  (2):\penalty0 1178--1203, 2024.

\bibitem[Shekhar et~al.(2022)Shekhar, Kim, and Ramdas]{ShekharKimRamdas2022}
Shubhanshu Shekhar, Ilmun Kim, and Aaditya Ramdas.
\newblock A permutation-free kernel two-sample test.
\newblock In \emph{Advances in Neural Information Processing Systems},
  volume~35, pages 18168--18180, 2022.

\bibitem[Shekhar et~al.(2023)Shekhar, Kim, and Ramdas]{ShekharKimRamdas2023}
Shubhanshu Shekhar, Ilmun Kim, and Aaditya Ramdas.
\newblock A permutation-free kernel independence test.
\newblock \emph{Journal of Machine Learning Research}, 24\penalty0
  (369):\penalty0 1--68, 2023.

\bibitem[Sz{\'e}kely and Rizzo(2013)]{SzekelyRizzo2013}
G{\'a}bor~J. Sz{\'e}kely and Maria~L. Rizzo.
\newblock Energy statistics: A class of statistics based on distances.
\newblock \emph{Journal of Statistical Planning and Inference}, 143\penalty0
  (8):\penalty0 1249--1272, 2013.

\bibitem[Sz{\'e}kely et~al.(2007)Sz{\'e}kely, Rizzo, and
  Bakirov]{SzekelyRizzoBakirov2007}
G{\'a}bor~J. Sz{\'e}kely, Maria~L. Rizzo, and Nail~K. Bakirov.
\newblock Measuring and testing dependence by correlation of distances.
\newblock \emph{The Annals of Statistics}, 35\penalty0 (6):\penalty0
  2769--2794, 2007.

\bibitem[Waudby-Smith et~al.(2024)Waudby-Smith, Arbour, Sinha, Kennedy, and
  Ramdas]{WaudbySmithArbourSinhaKennedyRamdas2024}
Ian Waudby-Smith, David Arbour, Ritwik Sinha, Edward~H. Kennedy, and Aaditya
  Ramdas.
\newblock Time-uniform central limit theory and asymptotic confidence
  sequences.
\newblock \emph{The Annals of Statistics}, 52\penalty0 (6):\penalty0
  2613--2640, 2024.

\bibitem[Whitehouse et~al.(2026)Whitehouse, Wu, and
  Ramdas]{WhitehouseWuRamdas2025}
Justin Whitehouse, Zhiwei~Steven Wu, and Aaditya Ramdas.
\newblock Time-uniform self-normalized concentration for vector-valued
  processes.
\newblock \emph{The Annals of Applied Probability}, 36\penalty0 (3):\penalty0
  1972--2013, 2026.

\bibitem[Zhang(2004)]{zhangASIP}
Li-xin Zhang.
\newblock Strong approximations of martingale vectors and their applications in
  {Markov}-chain adaptive designs.
\newblock \emph{Acta Mathematicae Applicatae Sinica, English Series},
  20\penalty0 (2):\penalty0 337--352, 2004.

\end{thebibliography}
\end{document}